# COVARIATE-ADJUSTED NONLINEAR REGRESSION[1]

By Xia Cui, Wensheng Guo,[2] Lu Lin and Lixing Zhu

*Shandong University, University of Pennsylvania, Shandong University,
Hong Kong Baptist University and East China Normal University*

In this paper, we propose a covariate-adjusted nonlinear regression model. In this model, both the response and predictors can only be observed after being distorted by some multiplicative factors. Because of nonlinearity, existing methods for the linear setting cannot be directly employed. To attack this problem, we propose estimating the distorting functions by nonparametrically regressing the predictors and response on the distorting covariate; then, nonlinear least squares estimators for the parameters are obtained using the estimated response and predictors. Root $n$-consistency and asymptotic normality are established. However, the limiting variance has a very complex structure with several unknown components, and confidence regions based on normal approximation are not efficient. Empirical likelihood-based confidence regions are proposed, and their accuracy is also verified due to its self-scale invariance. Furthermore, unlike the common results derived from the profile methods, even when plug-in estimates are used for the infinite-dimensional nuisance parameters (distorting functions), the limit of empirical likelihood ratio is still chi-squared distributed. This property eases the construction of the empirical likelihood-based confidence regions. A simulation study is carried out to assess the finite sample performance of the proposed estimators and confidence regions. We apply our method to study the relationship between glomerular filtration rate and serum creatinine.

Received March 2008; revised May 2008.

[1]Supported by a Grant (HKBU2030/07P) from Research Grants Council of Hong Kong, Hong Kong, China, NBRP (973 Program 2007CB814901) of China, NNSF project (10771123) of China, RFDP (20070422034) of China and NSF projects (Y2006A13 and Q2007A05) of Shandong Province of China.

[2]Supported in part by NCI R01 Grant CA84438.

**1. Introduction.** Consider the covariate-adjusted nonlinear regression model

$$(1) \qquad \begin{cases} Y = f(X, \boldsymbol{\beta}) + \varepsilon, \\ \widetilde{Y} = \psi(U)Y, \\ \widetilde{X}_r = \phi_r(U)X_r, \qquad r = 1, \dots, q, \end{cases}$$

where $Y$ is an unobservable response, $X = (X_1, \dots, X_q)^\tau$ is an unobservable predictor vector, $\boldsymbol{\beta}$ is the unknown $p \times 1$ vector parameter from a compact parameter space $\Theta \subset R^p$, $f(X, \boldsymbol{\beta})$ is a given continuous function in $X \in R^q$ and $\boldsymbol{\beta} \in R^p$, $\widetilde{Y}$ and $\widetilde{X}_r$ are the actual observable confounding variables and $\psi(U)$ and $\phi_r(U)$ are the unknown distorting functions of observable variable $U$. In this paper, we study point estimation and confidence region construction for the parameter $\boldsymbol{\beta}$. When $f(\cdot)$ is linear in $\boldsymbol{\beta}$, say, $f(X, \boldsymbol{\beta}) = \beta_0 + \sum_{l=1}^{p} \beta_l X_l$, Sentürk and Müller (2005) and Sentürk and Müller (2006), respectively, studied the consistency and asymptotic normality of the estimators and hypothesis testing. The model was motivated by an analysis of the regression of plasma fibrinogen concentration as response on serum transferrin level as predictor for 69 haemodialysis patients [see the details in Sentürk and Müller (2005)].

An example where a nonlinear model is relevant is the use of serum creatinine (SCr) to estimate glomerular filtration rate (GFR). GFR is traditionally considered the best overall index of renal function in health and disease. GFR can be measured by clearance techniques involving endogenous filtration markers (e.g., creatinine and urea) or exogenous markers (e.g., inulin, iohexol and iothalamate). Clearance studies of exogenous markers have been identified as the valid approach to measuring GFR across the spectrum of renal function. However, owing to the difficulty, complexity and expenses associated with measurement of GFR using clearance studies in clinical practice, clinicians traditionally estimate GFR from serum creatinine concentration. Because of the nonlinear relationship between GFR and SCr, numerous equations using different transformations of SCr have been developed, the most well-known being derived from the Modification of Diet in Renal Disease (MDRD) Study [see, e.g., Levey et al. (1999)]. It is well known that GFR and SCr are also related to body surface area [$\mathrm{BSA}(m^2) = 0.007184 * \mathrm{Kg}^{0.425} * \mathrm{cm}^{0.725}$]. Despite the nonlinear relationship between GFR and BSA, the adjusted GFR is usually calculated by simply dividing the unadjusted GFR over 1.73 times the BSA. The relationship between SCr and BSA is also nonlinear. Therefore, correcting the nonlinear effects of BSA on GFR and SCr may improve the estimated relationship between GFR and SCr, thereby leading to better estimation of GFR from SCr. As an illustration, we consider the data collected from the study A and study B of the MDRD study [see, e.g., Rosman et al. (1984) and Levey et al. (1994)]. Here, the relationship between GFR and SCr is of particular interest. We aim to show that the estimated relationship between GFR and SCr is more distinct after correcting for the distorting effect of BSA. Because the relationship of GFR and SCr is nonlinear, the covariate-adjusted linear regression model proposed by Sentürk and Müller (2005) cannot be used in our application. We will return to the example in Section 5.

To estimate the parameters in such nonlinear models, a natural consideration is whether the method proposed by Sentürk and Müller (2005) in the linear covariate-adjusted regression can be extended to the nonlinear setting. The basic idea of their method is to transform it into a varying-coefficient regression model, say,

$$(2) \qquad \widetilde{Y} = \alpha_0(U) + \sum_{l=1}^{p} \alpha_l(U)\widetilde{X}_l + \psi(U)\varepsilon,$$

where $\alpha_0(U) = \beta_0\psi(U)$, $\alpha_l(U) = \beta_l\frac{\psi(U)}{\phi_l(U)}$, $1 \le l \le p$. Note that $\beta_0 = \mathbb{E}[\alpha_0(U)]$ and $\beta_l = \frac{\mathbb{E}[\alpha_l(U)\widetilde{X}_l]}{\mathbb{E}[\widetilde{X}_l]}$, $1 \le l \le p$. Then, the bin method similar to that proposed in Fan and Zhang (2000) for longitudinal data is used to obtain consistent estimators $\hat{\alpha}_l$, $0 \le l \le p$. The estimators were further proved to be asymptotically normal by Sentürk and Müller (2006). The estimators of $\beta_0$ and $\beta_l$ can be obtained by replacing all the expectations by the corresponding weighted averages and $\alpha_l$ by $\hat{\alpha}_l$.

For model (1), when the nonlinear regression function $f(X, \boldsymbol{\beta})$ has a particular structure such as $f(X, \boldsymbol{\beta}) = f(\beta_0, \beta_1 X_1, \ldots, \beta_q X_q)$, we can let $\alpha_l(U) = \beta_l/\phi_l(U)$ and then transform the model to a nonlinear varying coefficient regression model to obtain estimators for $\alpha_l$. The relevant references are, among others [Staniswalis (2006) and Fan, Lin and Zhou (2006)]. The estimators for $\beta_l$ can be defined by using the fact $\beta_l = \frac{\mathbb{E}[\alpha_l(U)\widetilde{X}_l]}{\mathbb{E}[\widetilde{X}_l]}$. However, when $f$ is also related to $X_l^{a_l}$ for some given $a_l$, a similar approach cannot make $\beta_l$ estimable unless the expectations of $\phi_l(U)^{a_l}$ are given constants. When $a_l$ are also unknown parameters, we even have a problem to identify the parameters $\beta_l$ and $a_l$. Hence, it is obvious that, for more general models, we may not expect a useful, straightforward extension of this transformation method. This observation motivates us to develop a new approach to handle nonlinear models. In this paper, we consider a direct estimation procedure. We first use nonparametric regression to obtain consistent estimators of the distorting functions $\psi(\cdot)$ and $\phi_r(\cdot)$ by regressing the response and predictors on the distorting covariate, respectively, and obtain the estimates $(\hat{Y}, \hat{X}_1, \ldots, \hat{X}_q)$ for the unobservable response and predictors. We then apply least squares method to obtain the estimates of $\boldsymbol{\beta}$ in terms of $(\hat{Y}, \hat{X}_1, \ldots, \hat{X}_q)$. Note that the step of estimating the distorting functions is independent of the step of estimating $\boldsymbol{\beta}$. Thus, this method is obviously able to obtain a nonlinear least squares estimator of $\boldsymbol{\beta}$ for model (1).

There are two issues that need to be taken into account when the above procedure is applied. Although the above estimation procedure in principle is applicable, the root $n$-consistency and asymptotic normality are more difficult to study than the method suggested by Sentürk and Müller (2005) for linear models. This is due to the fact that nonlinear least squares estimator does not have a closed form and nonparametric estimators of the distorting functions have slower convergence rate than root $n$. We propose using the Delta method and $U$-statistics theory to obtain the asymptotic normality. The details are given in the Appendix. However, even though we can obtain the root $n$-consistency and asymptotic normality, the results in Section 2 show that the limiting variance of the estimator has a very complex structure; therefore, it is inconvenient to construct confidence region based on normal approximation. Hence, we propose using an empirical likelihood-based confidence region that avoids estimating the limiting variance and has better accuracy. Many advantages of empirical likelihood over the normal approximation-based method have been shown in the literature. In particular, it does not impose any prior constraints on the shape of the region, does not require the construction of a pivotal quantity and is range preserving and transformation respecting [see Hall and La Scala (1990)]. Owen (1991) applied the empirical likelihood to a linear regression model and proved that the empirical loglikelihood ratio is, asymptotically, a standard $\chi^2$-variable. Owen (2001) is a fairly comprehensive reference. Using empirical likelihood seems routine. However, a somewhat surprising result about our proposed method is that, although plug-in nonparametric estimators are used for infinite-dimensional nuisance distorting functions, the empirical likelihood ratio for the regression parameters is still of $\chi^2$ limit. This is very different from the existing literatures, because when nonparametric nuisance functions are replaced by plug-in estimators, the limit is often not tractable Chi-squares anymore, unless bias correction is implemented [see, e.g., Xue and Zhu (2007) and Zhu and Xue (2006) for details].

The rest of the paper is organized as follows. In Section 2, we describe the estimation procedure and the associated asymptotic results and discuss the efficiency of the estimators. In Section 3, we construct empirical likelihood based confidence regions for the parameters. In Section 4, some simulations are carried out to assess the performance of the proposed estimators and confidence regions. The application to the GFR and SCr data is presented in Section 5. Section 6 contains some concluding remarks. The technical proofs of all the asymptotic results are provided in the Appendix.

**2. Point estimation and asymptotic behavior.** As in Sentürk and Müller (2005), we assume that the mean distorting effect vanishes, that is,

(a)                $\mathbb{E}[\psi(U)] = 1,$        $\mathbb{E}[\phi_r(U)] = 1,$        $r = 1, \ldots, q.$

Besides, other basic assumptions are that

(b)     $(X_r, U, \varepsilon)$ are mutually independent,     $(Y, U)$ are independent,

(c)     $\mathbb{E}[\varepsilon] = 0,$     $\mathrm{Var}(\varepsilon) = \sigma^2.$

Assume that the available data are of the form $\{(U_i, \widetilde{X}_i, \widetilde{Y}_i), 1 \le i \le n\}$, for a sample of size $n$, where $\widetilde{X}_i = (\widetilde{X}_{1i}, \ldots, \widetilde{X}_{qi})^\tau$ are the $q$-dimensional observed predictors. The unobservable versions of $(\widetilde{X}_i, \widetilde{Y}_i)$ are denoted by $(X_i, Y_i)$. We write model (1) in a sample form $i = 1, \ldots, n$,

(3)
$$\begin{cases} Y_i = f(X_i, \boldsymbol{\beta}) + \varepsilon_i, \\ \widetilde{Y}_i = \psi(U_i) Y_i, \\ \widetilde{X}_{ri} = \phi_r(U_i) X_{ri}, \qquad r = 1, \ldots, q, \end{cases}$$

where $\{\varepsilon_i, 1 \le i \le n\}$ are independent and identically distributed random errors.

Now, our objective, based on the observations $\{(U_i, \widetilde{X}_i, \widetilde{Y}_i), 1 \le i \le n\}$, is to estimate the unknown parameter vector $\boldsymbol{\beta}$. From assumption (b),

(4)     $\psi(U) = \dfrac{\mathbb{E}[\widetilde{Y}|U]}{\mathbb{E}[Y]}, \qquad \phi_r(U) = \dfrac{\mathbb{E}[\widetilde{X}_r|U]}{\mathbb{E}[X_r]}, \qquad 1 \le r \le q.$

For convenience, we denote the density function of $U$ by $p(U)$ and define

(5)
$$g_Y(U) = \mathbb{E}[\widetilde{Y}|U] p(U),$$
$$g_r(U) = \mathbb{E}[\widetilde{X}_r|U] p(U), \qquad 1 \le r \le q.$$

There are some existing methods of estimating $\psi(U)$ and $\phi_r(U)$ [for details, see Eagleson and Müller (1997)]. Herein, we adopt a kernel method commonly used for ease of exposition. Then, for $1 \le r \le q$,

(6)
$$\hat{\psi}(u) = \frac{1/(nh) \sum_{i=1}^n K((u - U_i)/h) \widetilde{Y}_i}{1/(nh) \sum_{i=1}^n K((u - U_i)/h)} \times \frac{1}{\bar{\widetilde{Y}}} \triangleq \frac{\hat{g}_Y(u)}{\hat{p}(u)} \times \frac{1}{\bar{\widetilde{Y}}},$$
$$\hat{\phi}_r(u) = \frac{1/(nh) \sum_{i=1}^n K((u - U_i)/h) \widetilde{X}_{ri}}{1/(nh) \sum_{i=1}^n K((u - U_i)/h)} \times \frac{1}{\bar{\widetilde{X}}_r} \triangleq \frac{\hat{g}_r(u)}{\hat{p}(u)} \times \frac{1}{\bar{\widetilde{X}}_r},$$

where $\bar{\widetilde{Y}} = \frac{1}{n} \sum_{i=1}^n \widetilde{Y}_i$, $\bar{\widetilde{X}}_r = \frac{1}{n} \sum_{i=1}^n \widetilde{X}_{ri}$, $h$ is a bandwidth, and $K(\cdot)$ is a kernel function.

Let

(7)     $\hat{Y}_i = \widetilde{Y}_i / \hat{\psi}(U_i), \qquad \hat{X}_{ri} = \widetilde{X}_{ri} / \hat{\phi}_r(U_i) \quad \text{and} \quad \hat{X}_i = (\hat{X}_{1i}, \ldots, \hat{X}_{qi})^\tau.$

Then, the nonlinear least squares estimator $\hat{\boldsymbol{\beta}}$ for model (1) can be defined as the solution of the $p$ equations

(8)     $G_n^k(\boldsymbol{\beta}) = \sum_{i=1}^n (\hat{Y}_i - f(\hat{X}_i, \boldsymbol{\beta})) \dfrac{\partial f(\hat{X}_i, \boldsymbol{\beta})}{\partial \beta_k}$

for $k = 1, \ldots, p$, where $\partial f(\cdot, \boldsymbol{\beta})/\partial \beta_k$ is the partial derivative of $f$ with respect to $\beta_k$. Because $f$ is nonlinear in $\boldsymbol{\beta}$, no explicit solution can be calculated, and an iterative procedure is needed instead.

Observe that the asymptotic properties of $\hat{\boldsymbol{\beta}}$ defined by (8) depend on those of estimating functions $G_n^k(\boldsymbol{\beta})$. An approximation to $G_n^k(\boldsymbol{\beta})$ is provided in the following proposition, which is the key to obtaining the asymptotics of $\hat{\boldsymbol{\beta}}$ later on. The detailed proof is given in the Appendix.

PROPOSITION 1.   *Under conditions* (A1) *and* (A5)–(A8) *in the Appendix, we have, for* $1 \leq k \leq p$,

$$(9) \qquad n^{-1} G_n^k(\boldsymbol{\beta}) = n^{-1} R_n^k(\boldsymbol{\beta}) + o_P(n^{-1/2})$$

*and*

$$(10) \qquad n^{-1} G_n(\boldsymbol{\beta}) G_n^\tau(\boldsymbol{\beta}) = n^{-1} R_n(\boldsymbol{\beta}) R_n^\tau(\boldsymbol{\beta}) + o_P(1),$$

*where* $G_n(\boldsymbol{\beta}) = (G_n^1(\boldsymbol{\beta}), \ldots, G_n^p(\boldsymbol{\beta}))^\tau$, $R_n(\boldsymbol{\beta}) = (R_n^1(\boldsymbol{\beta}), \ldots, R_n^p(\boldsymbol{\beta}))^\tau$ *and*

$$
\begin{aligned}
R_n^k(\boldsymbol{\beta}) = {} & \sum_{i=1}^n \varepsilon_i f_{\beta_k}(X_i, \boldsymbol{\beta}) \\
& + \frac{1}{2} \sum_{i=1}^n \Bigg( (Y_i - \mathbb{E}[Y]) \frac{\mathbb{E}[Y f_{\beta_k}(X, \boldsymbol{\beta})]}{\mathbb{E}[Y]} \\
& \qquad\qquad - \sum_{l=1}^q (X_{li} - \mathbb{E}[X_l]) \frac{\mathbb{E}[X_l f_{x_l}(X, \boldsymbol{\beta}) f_{\beta_k}(X, \boldsymbol{\beta})]}{\mathbb{E}[X_l]} \Bigg) \\
& + \sum_{i=1}^n \Bigg( (\widetilde{Y}_i - Y_i) \frac{\mathbb{E}[Y f_{\beta_k}(X, \boldsymbol{\beta})]}{\mathbb{E}[Y]} \\
& \qquad\qquad - \sum_{l=1}^q (\widetilde{X}_{li} - X_{li}) \frac{\mathbb{E}[X_l f_{x_l}(X, \boldsymbol{\beta}) f_{\beta_k}(X, \boldsymbol{\beta})]}{\mathbb{E}[X_l]} \Bigg),
\end{aligned}
$$

*where* $f_{\beta_k}$, $f_{x_l}$ *denote the first derivative of* $f$ *with respect to* $\beta_k$ *and* $x_r$, *respectively,* $1 \leq k \leq p$, $1 \leq r \leq q$.

REMARK 1.   When specialized to the covariate-adjusted linear regression model, $f(X, \boldsymbol{\beta}) = \beta_0 + \sum_{l=1}^p \beta_l X_l$. The estimating equation $R_n^k(\boldsymbol{\beta})$ is

$$
\begin{aligned}
R_n^k(\boldsymbol{\beta}) = {} & \sum_{i=1}^n \varepsilon_i X_{ki} \\
& + \frac{1}{2} \sum_{i=1}^n \Bigg( (Y_i - \mathbb{E}[Y]) \frac{\mathbb{E}[Y X_k]}{\mathbb{E}[Y]} - \sum_{l=1}^q \beta_l (X_{li} - \mathbb{E}[X_l]) \frac{\mathbb{E}[X_k X_l]}{\mathbb{E}[X_l]} \Bigg) \\
& + \sum_{i=1}^n \Bigg( (\widetilde{Y}_i - Y_i) \frac{\mathbb{E}[Y X_k]}{\mathbb{E}[Y]} - \sum_{l=1}^q \beta_l (\widetilde{X}_{li} - X_{li}) \frac{\mathbb{E}[X_k X_l]}{\mathbb{E}[X_l]} \Bigg)
\end{aligned}
$$

and, for $\beta_0$, $X_0$ is equal to 1.

Root $n$-consistency and asymptotic normality of $\hat{\boldsymbol{\beta}}$ defined by $G_n^k(\boldsymbol{\beta})$ will be established in the following theorem, and the conditions required are given in the [Appendix](#). Here, we introduce the following notation ($1 \le s, k \le p$):

$$\Lambda(s,k) = \mathbb{E}[f_{\beta_r}(X, \boldsymbol{\beta}^0) f_{\beta_k}(X, \boldsymbol{\beta}^0)],$$

$$\zeta_k = \mathbb{E}[f_{\beta_k}(X, \boldsymbol{\beta}^0)], \qquad \zeta = (\zeta_1, \ldots, \zeta_p)^\tau,$$

$$\text{(11)} \qquad \eta_k = \frac{\mathbb{E}[Y f_{\beta_k}(X, \boldsymbol{\beta}^0)]}{\mathbb{E}[Y]}, \qquad \eta = (\eta_1, \ldots, \eta_p)^\tau,$$

$$\begin{aligned}
\Omega(s,k) = \mathbb{E}\Bigg\{ & \bigg( (Y - \mathbb{E}[Y]) \frac{\mathbb{E}[Y f_{\beta_s}(X, \boldsymbol{\beta}^0)]}{\mathbb{E}[Y]} \\
& - \sum_{l=1}^q (X_l - \mathbb{E}[X_l]) \frac{\mathbb{E}[X_l f_{x_l}(X, \boldsymbol{\beta}^0) f_{\beta_s}(X, \boldsymbol{\beta}^0)]}{\mathbb{E}[X_l]} \bigg) \\
& \times \bigg( (Y - \mathbb{E}[Y]) \frac{\mathbb{E}[Y f_{\beta_k}(X, \boldsymbol{\beta}^0)]}{\mathbb{E}[Y]} \\
& - \sum_{l=1}^q (X_l - \mathbb{E}[X_l]) \frac{\mathbb{E}[X_l f_{x_l}(X, \boldsymbol{\beta}^0) f_{\beta_k}(X, \boldsymbol{\beta}^0)]}{\mathbb{E}[X_l]} \bigg) \Bigg\},
\end{aligned}$$

$$\begin{aligned}
\Gamma(s,k) = \mathbb{E}\Bigg\{ & \bigg( (\widetilde{Y} - Y) \frac{\mathbb{E}[Y f_{\beta_s}(X, \boldsymbol{\beta}^0)]}{\mathbb{E}[Y]} \\
& - \sum_{l=1}^q (\widetilde{X}_l - X_l) \frac{\mathbb{E}[X_l f_{x_l}(X, \boldsymbol{\beta}^0) f_{\beta_s}(X, \boldsymbol{\beta}^0)]}{\mathbb{E}[X_l]} \bigg) \\
& \times \bigg( (\widetilde{Y} - Y) \frac{\mathbb{E}[Y f_{\beta_k}(X, \boldsymbol{\beta}^0)]}{\mathbb{E}[Y]} \\
& - \sum_{l=1}^q (\widetilde{X}_l - X_l) \frac{\mathbb{E}[X_l f_{x_l}(X, \boldsymbol{\beta}^0) f_{\beta_k}(X, \boldsymbol{\beta}^0)]}{\mathbb{E}[X_l]} \bigg) \Bigg\}.
\end{aligned}$$

THEOREM 1. *Let* $\hat{\boldsymbol{\beta}}$ *be defined by* [(8)](#). *If conditions* (A1)–(A8) *in the [Appendix](#) are satisfied, the following results hold:*

(i) $\hat{\boldsymbol{\beta}}$ *converges in probability to the true value* $\boldsymbol{\beta}^0$;

(ii)

$$\text{(12)} \qquad \sqrt{n}(\hat{\boldsymbol{\beta}} - \boldsymbol{\beta}^0) \xrightarrow{\mathcal{D}} N(0, \Sigma),$$

*where the covariance matrix can be represented as sum of $\Sigma = \sigma^2 \Lambda^{-1} + A + B + C$ with*

$$A = \tfrac{1}{4} \Lambda^{-1} \Omega \Lambda^{-1},$$
$$B = \Lambda^{-1} \Gamma \Lambda^{-1},$$
$$C = \sigma^2 \Lambda^{-1} (\tfrac{1}{2} \eta \zeta^\tau + \tfrac{1}{2} \zeta \eta^\tau) \Lambda^{-1}.$$

The proof is deferred to the [Appendix](). Let us consider the terms in the expression for the asymptotic covariance matrix in Theorem [1](). If there are no distortions with $\psi = \phi_r = 1$, and $Y = \widetilde{Y}$ and $X_r = \widetilde{X}_r$, then we can estimate $\boldsymbol{\beta}$ by least squares, minimizing $\sum_{i=1}^n (Y_i - f(X_i, \boldsymbol{\beta}))^2$ with respect to $\boldsymbol{\beta} \in \Theta$. The term $\sigma^2 \Lambda^{-1}$ equals the asymptotic covariance matrix of this least squares estimator. The matrix $A + B$ is caused by the distortions. The correlation between the first two terms of $R_n(\boldsymbol{\beta})$ leads to $C$.

REMARK 2. When the model is linear, the limiting variance matrix is reduced to

$$(13) \qquad \Sigma = \sigma^2 \Lambda^{-1} + A + B + \frac{\sigma^2}{2\mathbb{E}[Y]} (e_0 \boldsymbol{\beta}^{0\tau} + \boldsymbol{\beta}^0 e_0^\tau),$$

where $(0 \le s, k \le p$ and $X_0 = 1)$

$\Lambda(s, k) = \mathbb{E}[X_s X_k],$

$A(s, k) = \dfrac{1}{4} \beta_s^0 \beta_k^0 \mathbb{E} \left\{ \left( \dfrac{Y - \mathbb{E}[Y]}{\mathbb{E}[Y]} - \dfrac{X_s - \mathbb{E}[X_s]}{\mathbb{E}[X_s]} \right) \left( \dfrac{Y - \mathbb{E}[Y]}{\mathbb{E}[Y]} - \dfrac{X_k - \mathbb{E}[X_k]}{\mathbb{E}[X_k]} \right) \right\},$

$B(s, k) = \beta_s^0 \beta_k^0 \mathbb{E} \left\{ \left( \dfrac{\widetilde{Y} - Y}{\mathbb{E}[Y]} - \dfrac{\widetilde{X}_s - X_s}{\mathbb{E}[X_s]} \right) \left( \dfrac{\widetilde{Y} - Y}{\mathbb{E}[Y]} - \dfrac{\widetilde{X}_k - X_k}{\mathbb{E}[X_k]} \right) \right\},$

$e_0 = (1, 0, \ldots, 0)^\tau.$

Then, the asymptotic variance of $\sqrt{n}(\hat{\beta}_k - \beta_k^0)$ obtained from our new method is given by, $0 \le k \le p$

$$(V1): \begin{cases} \sigma_0^2 = \sigma^2 (\Lambda^{-1})_{00} + \dfrac{1}{4} \beta_0^{02} \mathbb{E} \left( \dfrac{Y - \mathbb{E}[Y]}{\mathbb{E}[Y]} \right)^2 + \beta_0^{02} \mathbb{E} \left( \dfrac{\widetilde{Y} - Y}{\mathbb{E}[Y]} \right)^2 + \dfrac{\sigma^2}{\mathbb{E}[Y]} \beta_0^0, \\ \sigma_k^2 = \sigma^2 (\Lambda^{-1})_{kk} + \dfrac{1}{4} \beta_k^{02} \mathbb{E} \left( \dfrac{Y - \mathbb{E}[Y]}{\mathbb{E}[Y]} - \dfrac{X_k - \mathbb{E}[X_k]}{\mathbb{E}[X_k]} \right)^2 \\ \qquad + \beta_k^{02} \mathbb{E} \left( \dfrac{\widetilde{Y} - Y}{\mathbb{E}[Y]} - \dfrac{\widetilde{X}_k - X_k}{\mathbb{E}[X_k]} \right)^2. \end{cases}$$

Now, we are in the position to make a comparison between the naive estimators $(\sigma_0^2, \sigma_k^2)$ and the transformation-based estimators $(\breve{\sigma}_0^2, \breve{\sigma}_k^2)$. According to the results of Sentürk and Müller (2006), the asymptotic variance of the transformation-based estimators is (with the correction of a typographical error in their paper)

$$(\text{V2}): \begin{cases} \breve{\sigma}_0^2 = \sigma^2 (\Lambda^{-1})_{00} + \sigma^2 (\Lambda^{-1})_{00} \operatorname{Var}[\psi(U)] + \beta_0^{02} \operatorname{Var}[\psi(U)], \\ \breve{\sigma}_k^2 = \sigma^2 (\Lambda^{-1})_{kk} + \sigma^2 (\Lambda^{-1})_{kk} \operatorname{Var}[\psi(U)] \\ \qquad + \beta_k^{02} \dfrac{\mathbb{E}[X_k^2]}{(\mathbb{E}[X_k])^2} \operatorname{Var}[\psi(U) - \phi_k(U)]. \end{cases}$$

The following result states a necessary and sufficient condition to judge when the naive estimators are more efficient with smaller limiting variance for this case. To formulate this result, we write $\Lambda$ as $\Lambda = (\Lambda_0, \Lambda_1, \ldots, \Lambda_p)$, with $\Lambda_k$ being $(p+1)$-dimensional column vector. Let $e_k = (0, \ldots, 0, 1, 0, \ldots, 0)^\tau$, which is of $p+1$ dimension and only $(k+1)$th element is 1.

THEOREM 2. *Treating as $\phi_0(\cdot) = 1$ and $X_0 = 1$, then $\sigma_k^2 \leq \breve{\sigma}_k^2$, if and only if,*

(14) $$\{\boldsymbol{\beta} : \boldsymbol{\beta}^\tau \mathbf{C}_k \boldsymbol{\beta} \leq 0, 0 \leq k \leq p\},$$

*where the matrix $\mathbf{C}_k$ is*

$$\mathbf{C}_k = \left\{ (-3 - 4\mathbb{E}[\psi^2(U)] + 8\mathbb{E}[\psi(U)\phi_k(U)]) \frac{\mathbb{E}[X_k^2]}{(\mathbb{E}[X_k])^2} \beta_k^2 \right.$$
$$\left. - 4\sigma^2 \operatorname{Var}[\psi(U)] (\Lambda^{-1})_{kk} \right\} \Lambda_0 \Lambda_0^\tau$$
$$+ \sigma^2 (-3 + 4\mathbb{E}[\psi^2(U)]) e_k e_k^\tau + 2\sigma^2 (e_k \Lambda_k^\tau + \Lambda_k e_k^\tau) I(k = 0)$$
$$+ (3 - 4\mathbb{E}[\psi(U)\phi_k(U)]) \frac{1}{\mathbb{E}[X_k]} \beta_k^2 (\Lambda_k \Lambda_0^\tau + \Lambda_0 \Lambda_k^\tau)$$
$$+ (-3 + 4\mathbb{E}[\psi^2(U)]) \beta_k^2 \Lambda.$$

*Furthermore, the matrices $\mathbf{C}_k$'s are symmetric with at least one negative eigenvalue.*

From this result, we realize that the limiting variance of the naive estimator relies on all the coefficients $\beta_k$; thus, no one can be uniformly more efficient than the other, in general. However, in the case without distortion, our method would perform worse than the transformation based method. The following examples show what the sufficient and necessary conditions for the naive estimator that has smaller limiting variance will reduce to, assuming that at least one variable in the covariate-adjusted linear model is contaminated.

Example 2.1.   Consider the simplest covariate-adjusted linear model with $p = 1$. Then,

$$(2.12) \qquad \begin{cases} Y = \beta_0 + X\beta_1 + \varepsilon, \\ \widetilde{Y} = \psi(U)Y, \\ \widetilde{X} = \phi(U)X. \end{cases}$$

The condition of Theorem 2 is

$$R = \{(\beta_0, \beta_1) : (\beta_0, \beta_1)\mathbf{C}_k(\beta_0, \beta_1)^\tau \le 0, 0 \le k \le 1\}$$

with the elements of matrices $\mathbf{C}_k$s defined by

$$\mathbf{C}_0(1,1) = \sigma^2\left(5 - 4\operatorname{Var}[\psi]\frac{(\mathbb{E}[X])^2}{\operatorname{Var}[X]}\right),$$

$$\mathbf{C}_0(1,2) = \mathbf{C}_0(2,1) = \sigma^2\mathbb{E}[X]\left(2 - 4\operatorname{Var}[\psi]\frac{\mathbb{E}[X^2]}{\operatorname{Var}[X]}\right),$$

$$\mathbf{C}_0(2,2) = \beta_0^2\operatorname{Var}[X](4\operatorname{Var}[\psi] + 1),$$

$$\mathbf{C}_1(1,1) = \beta_1^2(-3 - 4\mathbb{E}(\psi^2) + 8\mathbb{E}(\psi\phi))\frac{\operatorname{Var}[X]}{(\mathbb{E}[X])^2} - 4\sigma^2\operatorname{Var}[\psi]\frac{1}{\operatorname{Var}[X]},$$

$$\mathbf{C}_1(1,2) = \mathbf{C}_1(2,1) = -4\beta_1^2(\mathbb{E}[\psi^2] - \mathbb{E}[\psi\phi])\frac{\operatorname{Var}[X]}{\mathbb{E}[X]} - 4\sigma^2\operatorname{Var}[\psi]\frac{\mathbb{E}[X]}{\operatorname{Var}[X]},$$

$$\mathbf{C}_1(2,2) = \sigma^2\left(1 - 4\operatorname{Var}[\psi]\frac{(\mathbb{E}[X])^2}{\operatorname{Var}[X]}\right).$$

Unfortunately, the matrices $\mathbf{C}_k$s are rather complex and depend on the following many quantities: $\mathbb{E}[X]$, $\mathbb{E}[X^2]$, $\operatorname{Var}[\psi(U)]$, $\mathbb{E}[\psi(U)\phi(U)]$, $\sigma^2$ and the unknown parameters $\beta_0$ and $\beta_1$. It is not easy to directly imagine the region $R$. Figure 1 is helpful to understand the performance of two methods under the model assumptions that $\mathbb{E}[X] = 2$, $\mathbb{E}[X^2] = 5.144$, $\sigma^2 = 0.25$, $\operatorname{Var}[\psi(U)] = 0.08$ and $\mathbb{E}[\psi(U)\phi(U)] = 1$.

**3. Empirical likelihood-based confidence region.**   In this section, we discuss confidence region construction. Because we have obtained the asymptotic normality in the above section, normal approximation is a natural approach. However, as is shown in the theorem, the matrix $\Sigma$ is rather complex and includes several unknown components to be estimated. If we use normal approximation to construct a confidence region for $\boldsymbol{\beta}$, we need a plug-in estimation that involves the estimation for many unknown components. For comparison with the following empirical likelihood, we will still present normal approximation in the next section.

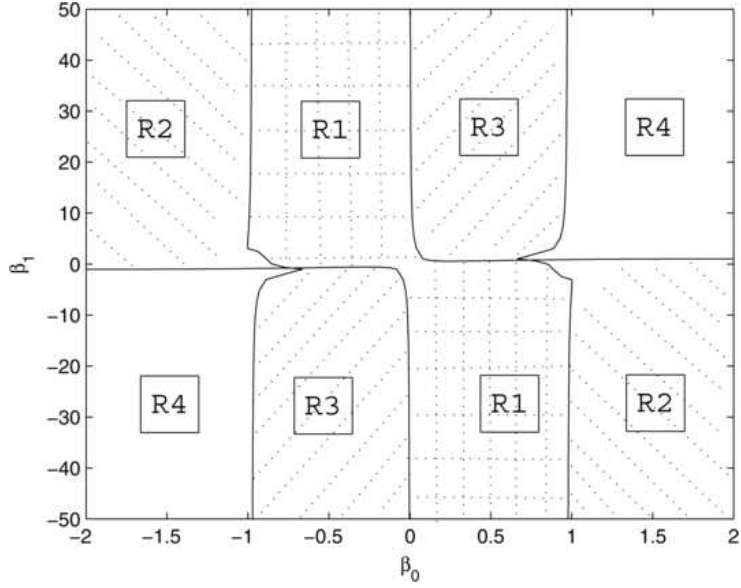

FIG. 1. *Region* R1 *corresponds to the case where the direct method performs better;* R2 *to that where the direct method performs better only for $\beta_1$;* R3 *to that where the direct method performs better only for $\beta_0$;* R4 *to that where Sentürk and Müller's method performs better.*

To construct a confidence region for $\beta$, we first introduce an auxiliary random variable $\xi_{n,i}^k(\boldsymbol{\beta}) = (Y_i - f(X_i, \boldsymbol{\beta})) \frac{\partial f(X_i, \boldsymbol{\beta})}{\partial \beta_k}$ and $\xi_{n,i}(\boldsymbol{\beta}) = (\xi_{n,i}^1(\boldsymbol{\beta}), \ldots, \xi_{n,i}^p(\boldsymbol{\beta}))^\tau$. Note that $\mathbb{E}[\xi_{n,i}(\boldsymbol{\beta})] = 0$ if $\boldsymbol{\beta} = \boldsymbol{\beta}^0$. Using this, an empirical log-likelihood ratio function is defined as

$$l_n(\boldsymbol{\beta}) = -2 \max \left\{ \sum_{i=1}^n \log(np_i) : p_i \geq 0, \sum_{i=1}^n p_i = 1, \sum_{i=1}^n p_i \xi_{n,i}(\boldsymbol{\beta}) = 0 \right\}.$$

Since $Y_i$ and $X_i$ in $l_n(\boldsymbol{\beta})$ are unobservable, a natural method is to replace them by the estimates $\hat{Y}_i$ and $\hat{X}_i$, respectively. A plug-in empirical log-likelihood, say $\hat{l}(\boldsymbol{\beta})$, is

$$(15) \quad \hat{l}(\boldsymbol{\beta}) = -2 \max \left\{ \sum_{i=1}^n \log(np_i) : p_i \geq 0, \sum_{i=1}^n p_i = 1, \sum_{i=1}^n p_i G_{n,i}(\boldsymbol{\beta}) = 0 \right\},$$

where $G_{n,i}(\boldsymbol{\beta}) = (G_{n,i}^1(\boldsymbol{\beta}), \ldots, G_{n,i}^p(\boldsymbol{\beta}))^\tau$ and $G_{n,i}^k(\boldsymbol{\beta}) = (\hat{Y}_i - f(\hat{X}_i, \boldsymbol{\beta})) \frac{\partial f(\hat{X}_i, \boldsymbol{\beta})}{\partial \beta_k}$, $k = 1, \ldots, p$. If we apply the standard method of Lagrange multipliers to find the optimal $p_i$ in $\hat{l}(\boldsymbol{\beta})$, the log empirical likelihood ratio is

$$(16) \qquad \hat{l}(\boldsymbol{\beta}) = 2 \sum_{i=1}^n \log\{1 + \lambda^\tau G_{n,i}(\boldsymbol{\beta})\},$$

where $\lambda$ is determined by

$$(17) \qquad \frac{1}{n} \sum_{i=1}^{n} \frac{G_{n,i}(\boldsymbol{\beta})}{1 + \lambda^\tau G_{n,i}(\boldsymbol{\beta})} = 0.$$

THEOREM 3. *Assume that conditions* (A1) *and* (A5)–(A8) *hold.* $\boldsymbol{\beta}^0$ *denotes the true parameter value. Then,*

$$(18) \qquad \hat{l}(\boldsymbol{\beta}^0) \xrightarrow{\mathcal{D}} \chi_p^2.$$

Based on the above theorem, an empirical likelihood confidence region for $\boldsymbol{\beta}$ with nominal confidence level $\alpha$ is

$$(19) \qquad I_{\alpha,\mathrm{EL}} = \{\boldsymbol{\beta} : \hat{l}(\boldsymbol{\beta}) \le c_\alpha\},$$

where $c_\alpha$ satisfies $P(\chi_p^2 \le c_\alpha) = 1 - \alpha$.

We note that, although the estimated $\hat{Y}_i$ and $\hat{X}_i$ that involve the estimation for infinite-dimensional nuisance parameters (distorting functions) are used, the Wilks' theorem still holds. This is very different from all of the existing results that need plug-in estimation for nuisance functions, unless bias correction is used. Xue and Zhu (2007) and Zhu and Xue (2006) show the necessity of bias correction in the relevant settings. This result is of importance for the use of empirical likelihood in semiparametric problems, because it tells us that we may not be able to give a unified theorem to deal with infinite-dimensional nuisance parameters. It merits further study.

We now consider normal approximation for constructing a confidence region for $\hat{\boldsymbol{\beta}}$. As we mentioned before, we have to estimate the unknowns in $\Sigma$ given in Theorem 1. Plugging the consistent estimators $\hat{\sigma}^2, \hat{\Lambda}, \hat{\zeta}, \hat{\eta}, \hat{\Gamma}$ and $\hat{\Omega}$ into $\Sigma$, we can obtain an estimator $\hat{\Sigma}$ of $\Sigma$. A confidence region with the approximate nominal coverage $1 - \alpha$ can be provided as

$$I_{\alpha,\mathrm{asy}} = \{\boldsymbol{\beta} : n(\hat{\boldsymbol{\beta}} - \boldsymbol{\beta})^\tau \hat{\Sigma}^{-1}(\hat{\boldsymbol{\beta}} - \boldsymbol{\beta}) \le c_\alpha\}.$$

**4. Simulation study.** In this section, we carry out simulations to investigate the performance of our proposed methods as outlined in Sections 2 and 3. Here, we choose a higher order kernel function $K(t) = \frac{15}{32}(3 - 7t^2)(1 - t^2)I(|t| \le 1)$ and use the leave-one-out cross-validation to select the optimal estimated bandwidth.

EXAMPLE 4.1. Consider the nonlinear regression model

$$(20) \qquad Y = \beta_1(1 - \exp(-\beta_2 X)) + \varepsilon,$$

where $\beta_1 = 4.309$, $\beta_2 = 0.208$, $\varepsilon \sim N(0, 2.0639)$ and $X$ is drawn from a normal distribution $\sim N(6.8, 26)$ truncated in interval $[0.45, 25]$. Assume that the distortion variable $U \sim N(0, 6)$ truncated in $[0, 6]$, and the distortion functions are $\psi(U) = (U + 1)^2/27.8170$ and $\phi(U) = (U + 1)/4.9459$, which satisfy the identity conditions. Then, 500 samples of $\widetilde{Y}$ and $\widetilde{X}$ were simulated from the specified distributions with sample sizes 200, 400 and 600. The estimated mean-squared errors for the estimators of $\beta_1$ and $\beta_2$ are (0.1983, 0.0717, 0.0315), (0.0022, 0.0006, 0.0006), respectively, for the sample sizes $n = 200, 400, 600$.

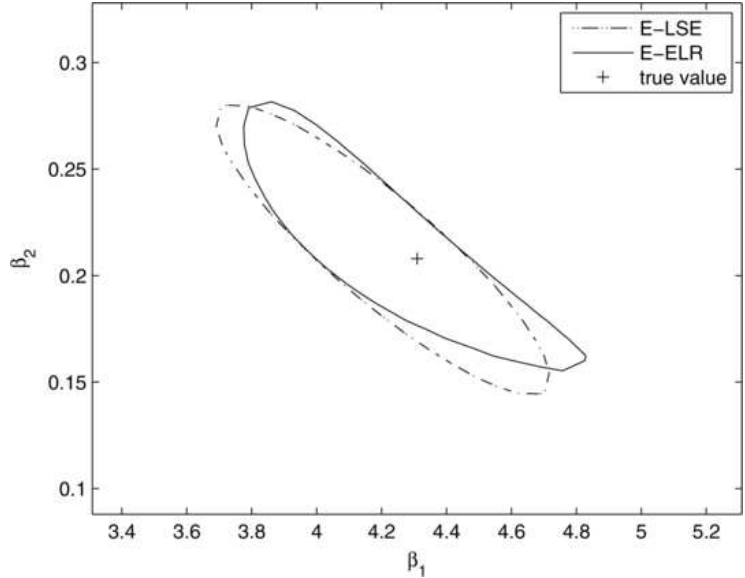

Fig. 2. *Confidence regions for Example 4.1. The solid and dash-dotted lines correspond to the empirical likelihood method and normal approximation method, respectively. The point "+" is the true value.*

Figure 2 reports a comparison of the two methods for constructing confidence regions, which are normal approximation and empirical likelihood, where the sample size is 400. The empirical likelihood provides a comparable confidence region as that of normal approximation. The coverage probabilities are reported in Table 1. We can see that the coverage probabilities of the empirical likelihood are uniformly higher than those of the normal approximation, which approach the nominal level as $n$ increases. Thus, in this example, the empirical likelihood is clearly superior to the normal approximation.

EXAMPLE 4.2. Assume that the underlying unobserved multiple regression model is

$$(21) \qquad Y = \beta_1(1 + X)^{\beta_2} + \varepsilon,$$

where $X$ arises from a uniform distribution with expectation 7/3 and variance 19/12, $\beta_1 = 2.5$, $\beta_2 = -1$ and $\varepsilon \sim N(0, 4)$. The confounding covariate $U$ was simulated from $\mathrm{Un}(4 - \sqrt{7}, 4 + \sqrt{7})$, and the distortion functions were chosen as $\psi(U) = (U + 10)^2/194.9160$ and $\phi(U) = (U + 34)/37.9160$, which satisfy the identifiability conditions. Again, we conducted 500 simulation runs with sample sizes 200, 400 and 600. The estimated mean-squared errors for the estimators of $\beta_1$ and $\beta_2$ are $(0.0121, 0.0054, 0.0038)$, $(0.0439, 0.0126, 0.0123)$, respectively, for the sample sizes $n = 200, 400$ and $600$. The confidence regions, constructed by two methods involving empirical likelihood method and normal approximation method, are shown in Figure 3 with the sample size being 400. We can see that the region based on the normal approximation is larger with the slightly lower coverage probabilities than those of the empirical likelihood, which are reported in Table 2. Thus, again, the empirical likelihood performs better.

5. **Application.** As an illustration of our method, we apply it to the baseline data collected from studies A and B of the Modification of Diet in Renal Disease (MDRD) Study [see Rosman et al. (1984) and Levey et al. (1994)]. The main goal of the original study was to demonstrate that dietary protein restriction can slow down the decline of the glomerular filtration rate (GFR). Here, we use the baseline unadjusted glomerular filtration rate (GFR) and serum creatinine (SCr) data to show that the relationship of the two variables can be substantially improved by correcting for the distorting effect of body surface area [BSA: $BSA(m^2) = 0.007184 * Kg^{0.425} * cm^{0.725}$]. After excluding a few missing data and couple outliers, the data set includes 819 subjects. Figure 4 depicts these original data. The main outcome in this example is the unadjusted GFR, which is measured as $mL/min$.

TABLE 1
*For Example 4.1, the coverage probabilities of the confidence regions on $(\beta_1, \beta_2)^\tau$*

| Method | 95% | | |
|---|---|---|---|
| | $n = 200$ | $n = 400$ | $n = 600$ |
| Empirical likelihood | 0.9140 | 0.9260 | 0.9360 |
| Normal approximation | 0.8520 | 0.8840 | 0.9340 |

TABLE 2
*For Example 4.2, the coverage probabilities of the confidence regions on $(\beta_1, \beta_2)^\tau$*

| Method | 95% | | |
|---|---|---|---|
| | $n = 200$ | $n = 400$ | $n = 600$ |
| Empirical likelihood | 0.9340 | 0.9360 | 0.9420 |
| Normal approximation | 0.9340 | 0.9260 | 0.9360 |

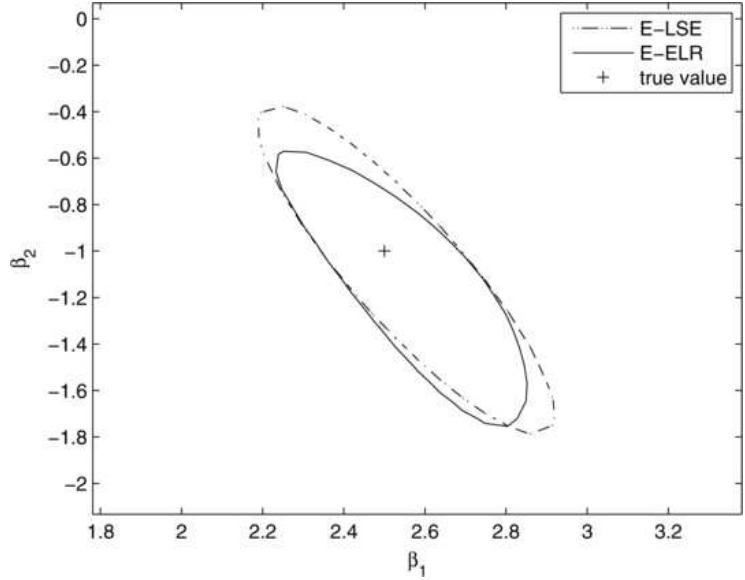

Fig. 3. *Confidence regions for Example 4.2. The solid and dash-dotted lines correspond to the empirical likelihood method and normal approximation method, respectively. The point "+" is the true value.*

For these data, let the observable response $\widetilde{Y}$ be the unadjusted GFR and the observable predictor $\widetilde{X}$ be the serum SCr. The consideration that the distortions of the response and predictor are unknown together with the fact that the relationship of GFR and SCr is nonlinear [see Levey et al. (1999) and Andrew et al. (2004)] motivates us to fit the serum creatinine data by a nonlinear model where the confounding variable $U$ is taken to be BSA. A relevant research [Ma et al. (2006)] showed that a nonlinear model for such kind of data is needed. Figure 5 presents the scatter plot of the logarithm of GFR against SCr, which shows a quadratic curve. Thus, we used a second order polynomial of SCr to fit an exponential model

$$(22) \qquad f(X, \boldsymbol{\beta}) = \beta_1 \exp(-\beta_2 X - \beta_3 X^2) + \beta_4.$$

A similar model, based on different transformations of SCr where they simply divided the unadjusted GFR by 1.73 times of BSA, was used in Levey et al. (1999). In this paper, instead of calculating the adjusted GFR by dividing it by 1.73 times the BSA, as is commonly done in medical literature, we allow the distorting function to be nonlinear and estimate it by kernel smoothing. The distorting function between SCr and BSA is also estimated using kernel smoothing. The kernel and bandwidth selector chosen are the same as those in Section 4. The increasing patterns of $\hat{\psi}(u)$ and $\hat{\phi}(u)$, which are displayed in Figures 6 and 7, provide evidence that $\psi(u)$ and $\phi(u)$ are not constant for these data. Therefore, the covariate-adjusted regression model is preferred over the existing methods that divide the GFR by 1.73 times the BSA. The restored data are then obtained by (7). We set the starting value of $\boldsymbol{\beta}$ as $(1, 0, 0, 1)^\tau$ and the final estimate is $\hat{\boldsymbol{\beta}} = (10.7100, 0.0759, -0.0004, 1.1528)^\tau$. The adjusted data using our proposed method and the nonlinear fitted curve are shown in Figure 8. After corrected for the BSA effects on both the response and predictor, the scatter plot shows that the adjusted data are closer to the fitted nonlinear curve, suggesting a more distinct relationship between GFR and SCr. This leads to a more accurate estimate of GFR using SCr.

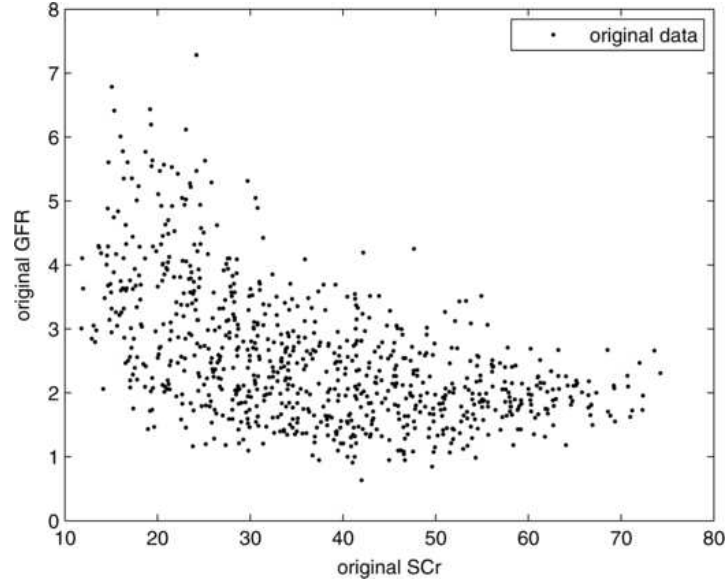

Fig. 4. *The scatter plot based on the data of original GFR and original SCr.*

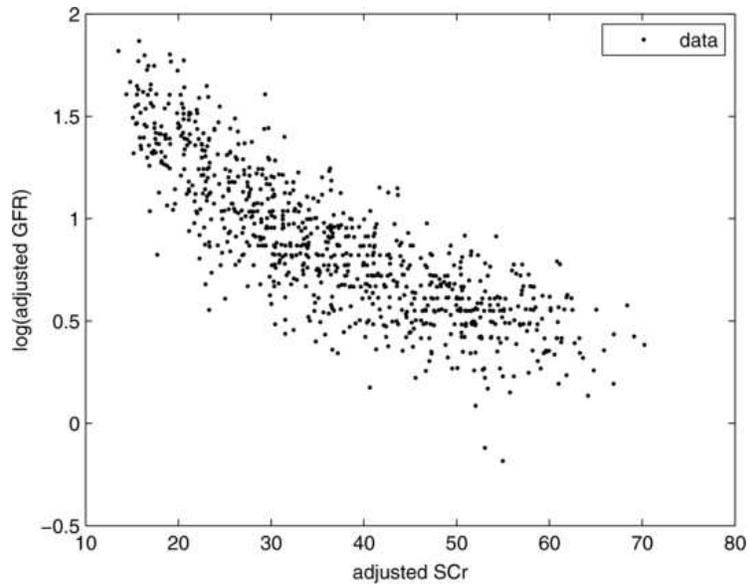

Fig. 5. *The scatter plot of logarithm of adjusted GFR against adjusted SCr.*

**6. Concluding remarks.**  For the models with distorting functions in response and predictors, we have proposed a direct approach to estimate and to make inference on the parameters of interest. For these models, there are several interesting issues that merit further studies. We recognize a merit in the original transformation method proposed by Sentürk and Müller (2005), which does not require divisions of the distorted variables by the estimated distorting functions, and the estimation may therefore be more numerically stable, compared to our proposed direct method, when some of the values of the distorting functions are closed to zero. How to inherit this merit in developing other more general estimation methods is an interesting topic. For the direct method, a practical remedy is to use their estimators plus a small constant, say, $1/n$. This is a commonly used remedy in nonparametric regression estimation [see, e.g., Zhu and Fang (1996)]. It will not damage their asymptotic properties. Another possible extension for the covariate-adjusted regression is to consider a model in which the distorted variables also have additive measurement errors. When the relationship between the underlying response and unobservable predictors is linear, it is possible to develop a method to estimate the parameters. However, the estimation may be a great challenge when the underlying relationship is nonlinear. Another related topic deserving further attention is on model checking when applying these methods to real data.

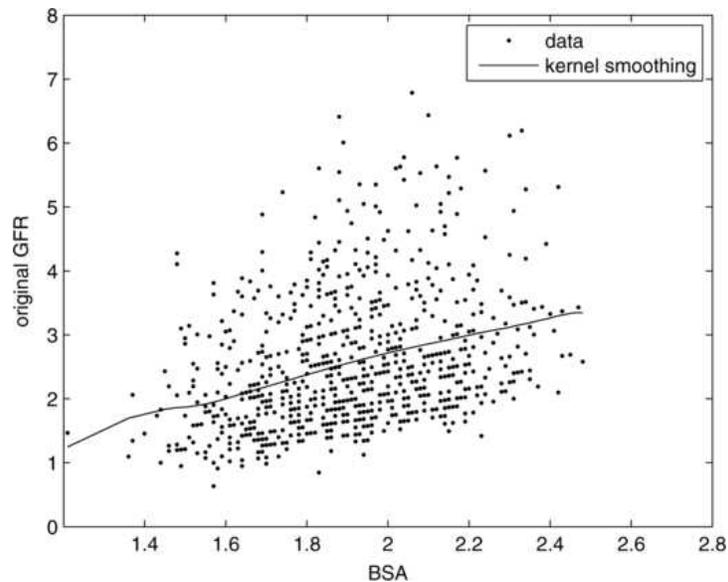

Fig. 6. *The kernel estimates based on the data of original GFR and BSA.*

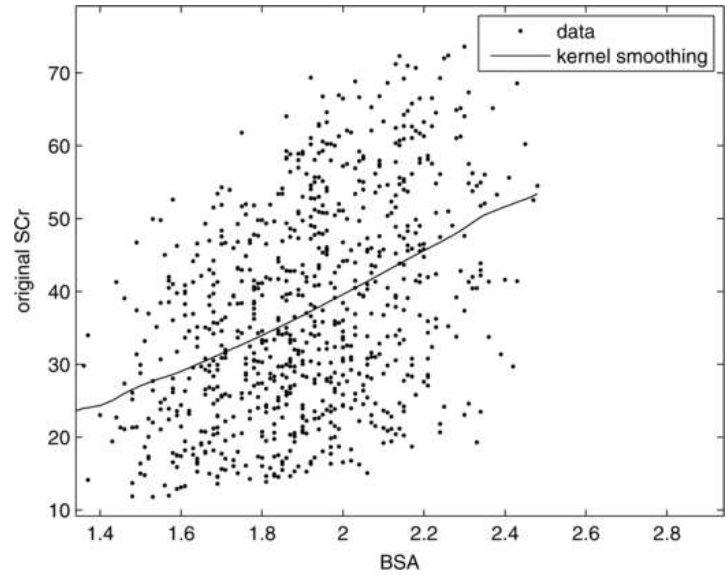

Fig. 7.   *The kernel estimates based on the data of original SCr and BSA.*

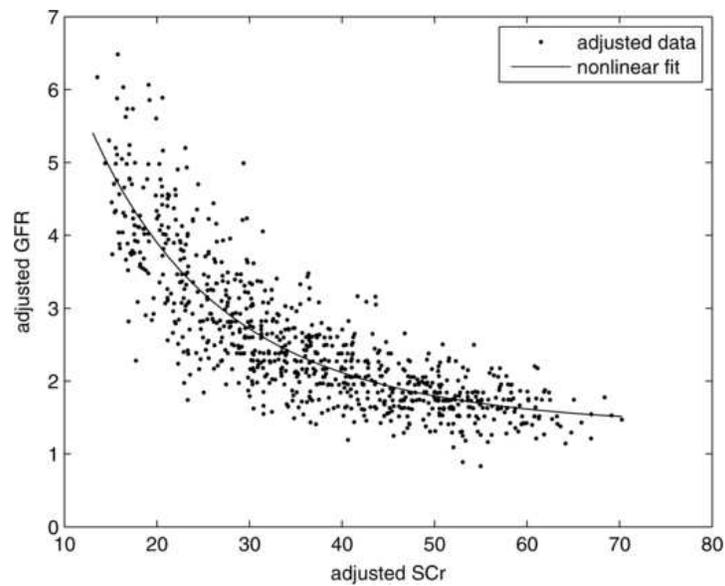

Fig. 8.   *The nonlinear fit based on model (22) and the data of adjusted GFR and adjusted SCr.*

## APPENDIX

The following conditions are needed for the results:

(A1) For all $s_1, s_2, s_3, s_4 = 0, 1, 2$, $s_1 + s_2 + s_3 + s_4 \leq 3$, $1 \leq k_1, k_2 \leq p$, $1 \leq l_1, l_2 \leq q$, the derivatives

$$\frac{\partial^{s_1+s_2+s_3+s_4} f(X, \boldsymbol{\beta})}{\partial^{s_1} \beta_{k_1} \partial^{s_2} \beta_{k_2} \partial^{s_3} x_{l_1} \partial^{s_4} x_{l_2}}$$

exist,

$$\left| \frac{\partial^{s_1+s_2+s_3+s_4} f(X, \boldsymbol{\beta})}{\partial^{s_1} \beta_{k_1} \partial^{s_2} \beta_{k_2} \partial^{s_3} x_{l_1} \partial^{s_4} x_{l_2}} \right| \leq C, \qquad \text{when } s_3 + s_4 \geq 1$$

and

$$E\left\{ \sup_{\boldsymbol{\beta}} \left| \frac{\partial^{s_1+s_2+s_3+s_4} f(X, \boldsymbol{\beta})}{\partial^{s_1} \beta_{k_1} \partial^{s_2} \beta_{k_2} \partial^{s_3} x_{l_1} \partial^{s_4} x_{l_2}} \right|^2 \right\} < \infty,$$

$$\text{when } 1 \leq s_1 + s_2 \leq 2 \text{ and } s_3 + s_4 = 0;$$

(A2) $\mathbb{E}[\varepsilon^4] < \infty$;

(A3) Suppose $\Lambda(s, k) = \mathbb{E}[f_{\beta_s}(X, \boldsymbol{\beta}^0) f_{\beta_k}(X, \boldsymbol{\beta}^0)] < \infty$ for $s, k = 1, \ldots, p$, and the matrix $\Lambda$ is positive definite;

(A4) $S(\boldsymbol{\beta}, \boldsymbol{\beta}_0) = \mathbb{E}[f(X, \boldsymbol{\beta}) - f(X, \boldsymbol{\beta}_0)]^2$ admits one unique minimum at $\boldsymbol{\beta} = \boldsymbol{\beta}_0$;

(A5) All $g_r(u) = \phi_r(u) p(u)$, $1 \leq r \leq q$, $g_Y(u) = \psi(u) p(u)$, $\phi_r(u)$, $\psi(u)$ and $p(u)$ are greater than a positive constant, are differential, and their derivatives satisfy the condition that there exits a neighborhood of the origin, say $\Delta$, and a constant $c > 0$ such that, for any $\delta \in \Delta$,

$$|g_r^{(3)}(u + \delta) - g_r^{(3)}(u)| \leq c|\delta|, \qquad 1 \leq r \leq q,$$

$$|g_Y^{(3)}(u + \delta) - g_Y^{(3)}(u)| \leq c|\delta|,$$

$$|p^{(3)}(u + \delta) - p^{(3)}(u)| \leq c|\delta|;$$

(A6) The continuous kernel function $K(\cdot)$ satisfies the following properties:

> (a1) the support of $K(\cdot)$ is the interval $[-1, 1]$,
> (a2) $K(\cdot)$ is symmetric about zero,
> (a3) $\int_{-1}^{1} K(u)\, du = 1$, $\int_{-1}^{1} u^i K(u)\, du = 0$, $i = 1, 2, 3$;

(A7) As $n \to \infty$, the bandwidth $h$ is in the range from $O(n^{-1/4} \log n)$ to $O(n^{-1/8})$;

(A8) $\mathbb{E}[Y]$ and $\mathbb{E}[X_r]$ are bounded away from 0, $\mathbb{E}[Y]^2 < \infty$ and $\mathbb{E}[X_r]^2 < \infty$.

These conditions are mild and can be satisfied in most of circumstances. Conditions (A1)–(A4) are essential for the asymptotic results of nonlinear least squares estimation. Condition (A5) is related to smoothness of the functions $g_r(\cdot)$ and $g_Y(\cdot)$ and the density function $p(\cdot)$ of $U$. Conditions (A6)–(A7) are commonly assumed for the root $n$ consistency of kernel based estimation [see, e.g., Zhu and Fang (1996) and Stoker and Härdle (1989)]. Condition (A8) is special for this problem see Şentürk and Müller (2006)].

We start with a lemma, which is frequently used in the process of the proof. Then, we proceed to the proof of Proposition 1 and Theorem 1, concerning the consistency and asymptotic normality.

LEMMA A.1. *Let $\eta(x)$ be a continuous function satisfying $\mathbb{E}[\eta(X)]^2 < \infty$. Assume that conditions* (A5)–(A8) *hold. The following asymptotic representation holds:*

$$\sum_{i=1}^{n} (\hat{Y}_i - Y_i)\eta(X_i) = \frac{1}{2} \sum_{i=1}^{n} (Y_i - \mathbb{E}[Y]) \frac{\mathbb{E}[Y\eta(X)]}{\mathbb{E}[Y]}$$

$$+ \sum_{i=1}^{n} (\widetilde{Y}_i - Y_i) \frac{\mathbb{E}[Y\eta(X)]}{\mathbb{E}[Y]} + o_P(\sqrt{n}).$$

*When replacing $\hat{Y}_i$, $Y_i$, $\widetilde{Y}_i$ with $\hat{X}_{ri}$, $X_{ri}$ and $\widetilde{X}_{ri}$, respectively, we can have similar asymptotic representation.*

PROOF OF LEMMA A.1. Recall that $\hat{Y}_i = \widetilde{Y}_i / \hat{\psi}(U_i)$, which is given by (7) and (6). Simple calculations give the expression

$$\sum_{i=1}^{n} (\hat{Y}_i - Y_i)\eta(X_i) = \sum_{i=1}^{n} \hat{Y}_i \eta(X_i) - \sum_{i=1}^{n} Y_i \eta(X_i).$$

Thus, it suffices to deal with the quantity $\sum_{i=1}^{n} \hat{Y}_i \eta(X_i)$. According to (6), it holds that

$$\sum_{i=1}^{n} \hat{Y}_i \eta(X_i)$$

(23)
$$= \sum_{i=1}^{n} \widetilde{Y}_i \eta(X_i) \left( \frac{\hat{p}(U_i)}{\hat{g}_Y(U_i)} \right) \bar{\tilde{Y}}$$

$$= \sum_{i=1}^{n} \widetilde{Y}_i \eta(X_i) \left( \frac{\hat{p}(U_i)}{\hat{g}_Y(U_i)} \right) \mathbb{E}[Y] + \sum_{i=1}^{n} \widetilde{Y}_i \eta(X_i) \left( \frac{\hat{p}(U_i)}{\hat{g}_Y(U_i)} \right) (\bar{\tilde{Y}} - \mathbb{E}[Y])$$

$$\triangleq L_1 + L_2.$$

First, we analyze $L_1$. The proof is divided into three steps. Denote by $L_{11}$, $L_{12}$ and $L_{13}$ the quantities

(24)
$$L_{11} = \sum_{i=1}^{n} Y_i \eta(X_i), \qquad L_{12} = \sum_{i=1}^{n} Y_i \eta(X_i) \frac{\hat{g}_Y(U_i)}{g_Y(U_i)},$$

$$L_{13} = \sum_{j=1}^{n} Y_i \eta(X_i) \frac{\hat{p}(U_i)}{p(U_i)}.$$

STEP 1. Show that

(25)
$$L_1 = L_{11} - L_{12} + L_{13} + o_P(\sqrt{n}).$$

Applying the equation $\psi(\cdot) = \frac{g_Y(\cdot)}{p(\cdot)EY}$, we have

$$L_1 = \sum_{i=1}^{n} Y_i \eta(X_i) \frac{g_Y(U_i)}{p(U_i)} \frac{\hat{p}(U_i)}{\hat{g}_Y(U_i)}$$

$$= \sum_{i=1}^{n} Y_i \eta(X_i) \frac{g_Y(U_i)}{p(U_i)} \frac{p(U_i)}{g_Y(U_i)} \left( 1 - \frac{\hat{g}_Y(U_i)}{g_Y(U_i)} + \frac{\hat{p}(U_i)}{p(U_i)} \right)$$

$$\quad + L_1^{R_1} - L_1^{R_2}$$

$$= \sum_{i=1}^{n} Y_i \eta(X_i) - \sum_{i=1}^{n} Y_i \eta(X_i) \frac{\hat{g}_Y(U_i)}{g_Y(U_i)} + \sum_{i=1}^{n} Y_i \eta(X_i) \frac{\hat{p}(U_i)}{p(U_i)}$$

$$\quad + L_1^{R_1} - L_1^{R_2}$$

$$= L_{11} - L_{12} + L_{13} + L_1^{R_1} - L_1^{R_2},$$

where

$$L_1^{R_1} = \sum_{i=1}^{n} Y_i \eta(X_i) \frac{(\hat{g}_Y(U_i) - g_Y(U_i))^2}{g_Y(U_i)\hat{g}_Y(U_i)},$$

$$L_1^{R_2} = \sum_{i=1}^{n} Y_i \eta(X_i) \frac{(\hat{g}_Y(U_i) - g_Y(U_i))(\hat{p}(U_i) - p(U_i))}{p(U_i)\hat{g}_Y(U_i)}.$$

Equation (25) can be concluded by proving

$$L_1^{R_i} = o_P(n^{1/2}), \qquad 1 \le i \le 2. \tag{26}$$

Invoking Theorem 2.1.8 of Rao (1983), Lemma 3 of Zhu and Fang (1996), condition (A7) and the Law of Large Numbers (LLN) for $\frac{1}{n}\sum_{i=1}^n Y_i\eta(X_i)$, we then follow the arguments used in Zhu and Fang (1996) to yield the higher order term

$$L_1^{R_1} = O_P\left((h^4 + n^{-1/2}h^{-1}\log n)^2\left\{\sum_{i=1}^n Y_i\eta(X_i)\right\}\right)$$

$$= o_P(n^{1/2}).$$

Similar arguments yield (26) for $i = 2$.

STEP 2. Show that

$$L_{12} = \frac{1}{2}\sum_{i=1}^n Y_i\eta(X_i) + \frac{1}{2}\sum_{i=1}^n Y_i\frac{\mathbb{E}[Y\eta(X)]}{\mathbb{E}[Y]} + o_P(n^{1/2}),$$

$$\tag{27}$$

$$L_{13} = \frac{1}{2}\sum_{i=1}^n Y_i\eta(X_i) + \frac{1}{2}\sum_{i=1}^n \mathbb{E}[Y\eta(X)] + o_P(n^{1/2}).$$

To analyze $L_{12}$, we need to work on

$$\sum_{i=1}^n Y_i\eta(X_i)\frac{\hat{g}_Y(U_i)}{g_Y(U_i)}.$$

Again, applying the similar arguments used by Zhu and Fang (1996), we approximate the above sum by a $U$-statistic with a varying kernel with the bandwidth $h$. We can then have the asymptotic representation [see, e.g., Serfling (1980)] as

$$\sum_{i=1}^n Y_i\eta(X_i)\frac{\hat{g}_Y(U_i)}{g_Y(U_i)}$$

$$= \frac{1}{2}\sum_{i=1}^n \frac{1}{h}\int Y_i\eta(X_i)\frac{1}{g_Y(U_i)}K\left(\frac{U_i-u}{h}\right)y\psi(u)p_{Y,U}(y,u)\,dy\,du$$

$$+ \frac{1}{2}\sum_{j=1}^n \frac{1}{h}\int Y_i\psi(U_i)\frac{1}{g_Y(u)}K\left(\frac{U_i-u}{h}\right)tp_{T,U}(t,u)\,dt\,du + o_P(n^{1/2})$$

with $p_{Y,U}(y,u)$ being the density function of $(Y,U)$, $T = Y\eta(X)$ and $p_{T,U}(t,u)$ the density of variable $(T,U)$. Note that $Y$ and $U$ are independent, $T$ and

$U$ are independent. Invoking conditions (A5)–(A8), we obtain

$$\sum_{i=1}^{n} Y_i \eta(X_i) \frac{\hat{g}_Y(U_i)}{g_Y(U_i)}$$

$$= \frac{1}{2} \sum_{i=1}^{n} Y_i \eta(X_i) \frac{\mathbb{E}[Y]}{g_Y(U_i)} \frac{1}{h} \int K\left(\frac{U_i - u}{h}\right) \psi(u) p_U(u)\, du$$

$$+ \frac{1}{2} \sum_{i=1}^{n} Y_i \psi(U_i) \mathbb{E}[Y \eta(X)] \frac{1}{h} \int \frac{p_U(u)}{g_Y(u)} K\left(\frac{U_i - u}{h}\right) du + o_P(n^{1/2})$$

$$= \frac{1}{2} \sum_{i=1}^{n} Y_i \eta(X_i) + \frac{1}{2} \sum_{i=1}^{n} Y_i \frac{\mathbb{E}[Y \eta(X)]}{\mathbb{E}[Y]} + o_P(n^{1/2}).$$

The result of (27) for $L_{13}$ can be similarly proved. Combining (25) with (27), we have

$$
\begin{aligned}
(28) \quad L_1 &= \sum_{i=1}^{n} (Y_i \eta(X_i) - \mathbb{E}[Y \eta(X)]) \\
&\quad - \frac{1}{2} \sum_{i=1}^{n} (Y_i - \mathbb{E}[Y]) \frac{\mathbb{E}[Y \eta(X)]}{\mathbb{E}[Y]} \\
&\quad + n\mathbb{E}[Y \eta(X)] + o_P(n^{1/2}).
\end{aligned}
$$

STEP 3.  Show that

$$(29) \qquad L_2 = \sum_{i=1}^{n} (\widetilde{Y}_i - \mathbb{E}[Y]) \frac{\mathbb{E}[Y \eta(X)]}{\mathbb{E}[Y]} + o_P(n^{1/2}).$$

From (28) and the definition of $L_2$ in (23), we derive that,

$$
\begin{aligned}
L_2 &= \frac{\bar{\tilde{Y}} - \mathbb{E}[Y]}{\mathbb{E}[Y]} L_1 \\
&= \sum_{i=1}^{n} (Y_i \eta(X_i) - \mathbb{E}[Y \eta(X)]) \frac{\bar{\tilde{Y}} - \mathbb{E}[Y]}{\mathbb{E}[Y]} \\
&\quad - \frac{1}{2} \sum_{i=1}^{n} (Y_i - \mathbb{E}[Y]) \frac{\mathbb{E}[Y \eta(X)]}{\mathbb{E}[Y]} \frac{\bar{\tilde{Y}} - \mathbb{E}[Y]}{\mathbb{E}[Y]} \\
&\quad + n\mathbb{E}[Y \eta(X)] \frac{\bar{\tilde{Y}} - \mathbb{E}[Y]}{\mathbb{E}[Y]} + o_P(n^{1/2}) \\
&= \sum_{i=1}^{n} (\widetilde{Y}_i - \mathbb{E}[Y]) \frac{\mathbb{E}[Y \eta(X)]}{\mathbb{E}[Y]} + o_P(n^{1/2}),
\end{aligned}
$$

where the last equality is obtained by applying LLN to $\frac{\bar{\tilde{Y}}-\mathbb{E}[Y]}{\mathbb{E}[Y]}$, $\frac{1}{n}\sum_{i=1}^{n}(Y_i\eta(X_i)-\mathbb{E}[Y\eta(X)])$ and $\frac{1}{n}\sum_{i=1}^{n}(Y_i-\mathbb{E}[Y])$.

Finally, together with the asymptotic representation of $L_1$ and $L_2$ in (28) and (29), the desired result is easy to arrive at.  □

PROOF OF PROPOSITION 1.

PROOF OF (9).   Start by writing $G_n^k(\boldsymbol{\beta}) = D_1 + D_2 + D_3$ by using Taylor expansion, where

$$D_1 = \sum_{i=1}^{n}\varepsilon_i\frac{\partial f(X_i,\boldsymbol{\beta})}{\partial\beta_k} + \sum_{i=1}^{n}(\hat{Y}_i-Y_i)\frac{\partial f(X_i,\boldsymbol{\beta})}{\partial\beta_k}$$
$$- \sum_{i=1}^{n}\sum_{l=1}^{q}(\hat{X}_{li}-X_{li})\frac{\partial f(X_i,\boldsymbol{\beta})}{\partial x_l}\frac{\partial f(X_i,\boldsymbol{\beta})}{\partial\beta_k},$$

$$D_2 = -\frac{1}{2}\sum_{i=1}^{n}\left(\sum_{l=1}^{q}\sum_{r=1}^{q}\frac{\partial^2 f(X_i^*,\boldsymbol{\beta})}{\partial x_l\,\partial x_r}(\hat{X}_{li}-X_{li})(\hat{X}_{ri}-X_{ri})\right)$$
$$\times\left(\frac{\partial f(X_i,\boldsymbol{\beta})}{\partial\beta_k} + \sum_{l=1}^{q}\frac{\partial^2 f(X_i^{**},\boldsymbol{\beta})}{\partial\beta_k\,\partial x_l}(\hat{X}_{li}-X_{li})\right),$$

$$D_3 = \sum_{i=1}^{n}\left(\varepsilon_i + (\hat{Y}_i-Y_i) - \sum_{l=1}^{q}\frac{\partial f(X_i,\boldsymbol{\beta})}{\partial x_l}(\hat{X}_{li}-X_{li})\right)$$
$$\times\left(\sum_{l=1}^{q}\frac{\partial^2 f(X_i^{**},\boldsymbol{\beta})}{\partial\beta_k\,\partial x_l}(\hat{X}_{li}-X_{li})\right),$$

where $X_i^{**} = (X_{1i}^{**},\ldots,X_{qi}^{**})^{\tau}$ and $X_i^* = (X_{1i}^*,\ldots,X_{qi}^*)^{\tau}$ with both $X_{li}^{**}$ and $X_{li}^*$ are two points between $\hat{X}_{li}$ and $X_{li}$. Applying Lemma A.1 to the second term of $D_1$ with $\eta(X) = \frac{\partial f(X,\boldsymbol{\beta})}{\partial\beta_k}$, and to the third term with $\eta(X) = \frac{\partial f(X,\boldsymbol{\beta})}{\partial x_l}\frac{\partial f(X,\boldsymbol{\beta})}{\partial\beta_k}$, we can conclude that

$$\sum_{i=1}^{n}(\hat{Y}_i-Y_i)\frac{\partial f(X_i,\boldsymbol{\beta})}{\partial\beta_k}$$
$$= \frac{1}{2}\sum_{i=1}^{n}(Y_i-\mathbb{E}[Y])\frac{\mathbb{E}[Y(\partial f(X,\boldsymbol{\beta}))/\partial\beta_k]}{\mathbb{E}[Y]}$$
$$+ \sum_{i=1}^{n}(\tilde{Y}_i-Y_i)\frac{\mathbb{E}[Y(\partial f(X,\boldsymbol{\beta}))/\partial\beta_k]}{\mathbb{E}[Y]} + o_P(\sqrt{n})$$

$$= \frac{1}{2}\sum_{i=1}^{n}(Y_i - \mathbb{E}[Y])\frac{\mathbb{E}[Y(\partial f(X,\boldsymbol{\beta}))/\partial\beta_k]}{\mathbb{E}[Y]}$$

$$+ \sum_{i=1}^{n}(\widetilde{Y}_i - Y_i)\frac{\mathbb{E}[Y(\partial f(X,\boldsymbol{\beta}))/\partial\beta_k]}{\mathbb{E}[Y]} + o_P(\sqrt{n}),$$

$$\sum_{i=1}^{n}\sum_{l=1}^{q}(\hat{X}_{li} - X_{li})\frac{\partial f(X_i,\boldsymbol{\beta})}{\partial x_l}\frac{\partial f(X_i,\boldsymbol{\beta})}{\partial\beta_k}$$

$$= \frac{1}{2}\sum_{i=1}^{n}\sum_{l=1}^{q}(X_{li} - \mathbb{E}[X_l])\frac{\mathbb{E}[X_l(\partial f(X,\boldsymbol{\beta}))/\partial\beta_k(\partial f(X,\boldsymbol{\beta}))/\partial x_l]}{\mathbb{E}[X_l]}$$

$$+ \sum_{i=1}^{n}\sum_{l=1}^{q}(\widetilde{X}_{li} - X_{li})\frac{\mathbb{E}[X_l(\partial f(X,\boldsymbol{\beta}))/\partial\beta_k(\partial f(X,\boldsymbol{\beta}))/\partial x_l]}{\mathbb{E}[X_l]} + o_P(\sqrt{n}).$$

To accomplish the proof, we only need to show that both $D_2$ and $D_3$ are asymptotically negligible. Note the arguments, which can be found in Zhu and Fang (1996),

$$\sup_u |\hat{p}(u) - p(u)| = O(h^4 + n^{-1/2}h^{-1}\log n), \qquad \text{a.s.},$$

(30)

$$\sup_u |\hat{g}_Y(u) - g_Y(u)| = O_P(h^4 + n^{-1/2}h^{-1}\log n).$$

Combining this with conditions (A1) and (A7)–(A8), we have

$$D_2 = O_P\left((h^4 + n^{-1/2}h^{-1}\log n)^2\left\{\sum_{i=1}^{n}\sum_{l=1}^{q}\sum_{r=1}^{q}C\left|\frac{\partial f(X_i,\boldsymbol{\beta})}{\partial\beta_k}X_{ri}X_{li}\right|\right\}\right)$$

$$+ O_P\left((h^4 + n^{-1/2}h^{-1}\log n)^3\left\{\sum_{i=1}^{n}\sum_{l=1}^{q}\sum_{r=1}^{q}C\right\}\right)$$

$$= O_P((h^4 + n^{-1/2}h^{-1}\log n)^2 n) + O_P((h^4 + n^{-1/2}h^{-1}\log n)^3 n)$$

$$= o_P(n^{1/2}).$$

Similarly, $D_3$ is $o_P(n^{1/2})$. The desired result follows from above analysis.

PROOF OF (10). Note that $R_n(\boldsymbol{\beta}) = O_P(n^{1/2})$ and

$$\frac{1}{n}G_n(\boldsymbol{\beta})G_n^\tau(\boldsymbol{\beta})$$

$$= \frac{1}{n}R_n(\boldsymbol{\beta})R_n^\tau(\boldsymbol{\beta}) + \frac{1}{n}(G_n(\boldsymbol{\beta}) - R_n(\boldsymbol{\beta}))R_n^\tau(\boldsymbol{\beta})$$

$$+ \frac{1}{n}R_n(\boldsymbol{\beta})(G_n(\boldsymbol{\beta}) - R_n(\boldsymbol{\beta}))^\tau + \frac{1}{n}(G_n(\boldsymbol{\beta}) - R_n(\boldsymbol{\beta}))(G_n(\boldsymbol{\beta}) - R_n(\boldsymbol{\beta}))^\tau.$$

Combining this with (9), it can immediately lead to the desired result of (10). □

Proof of Theorem 1.

Proof of (i). Denote by $S_n(\boldsymbol{\beta})$ the quantity $S_n(\boldsymbol{\beta}|\{X_i, Y_i\}) = \sum_{i=1}^n (Y_i - f(X_i, \boldsymbol{\beta}))^2$. The main point of the proof lies in stating that, for any $\boldsymbol{\beta}$ in $\Theta$, $S_n(\boldsymbol{\beta}|\{\hat{X}_i, \hat{Y}_i\}) - S_n(\boldsymbol{\beta}|\{X_i, Y_i\}) = o_P(n)$. This implies, by the fact that $\hat{\boldsymbol{\beta}} = \arg\min_{\boldsymbol{\beta} \in \Theta} S_n(\boldsymbol{\beta}|\{\hat{X}_i, \hat{Y}_i\})$ and the arguments used in Wu (1981) under conditions (A1)–(A4), that $S_n(\boldsymbol{\beta}|\{X_i, Y_i\})$ converges in probability to $S(\boldsymbol{\beta}, \boldsymbol{\beta}^0)$, admitting a unique minimum at $\boldsymbol{\beta} = \boldsymbol{\beta}^0$. The consistency of $\hat{\boldsymbol{\beta}}$ follows from Lemma 1 of Wu (1981). Now, after simple calculations, we can obtain the decomposition $S_n(\boldsymbol{\beta}|\{\hat{X}_i, \hat{Y}_i\}) - S_n(\boldsymbol{\beta}|\{X_i, Y_i\}) = F_1 + F_2 + F_3 + F_4 + F_5$, where

$$F_1 = \sum_{i=1}^n (\hat{Y}_i - Y_i)^2,$$

$$F_2 = \sum_{i=1}^n (f(\hat{X}_i, \boldsymbol{\beta}) - f(X_i, \boldsymbol{\beta}))^2,$$

$$F_3 = 2\sum_{i=1}^n (\hat{Y}_i - Y_i)\varepsilon_i,$$

$$F_4 = -2\sum_{i=1}^n (f(\hat{X}_i, \boldsymbol{\beta}) - f(X_i, \boldsymbol{\beta}))(\hat{Y}_i - Y_i),$$

$$F_5 = -2\sum_{i=1}^n (f(\hat{X}_i, \boldsymbol{\beta}) - f(X_i, \boldsymbol{\beta}))\varepsilon_i.$$

Note that $f(\hat{X}_i, \boldsymbol{\beta}) - f(X_i, \boldsymbol{\beta}) = \sum_{l=1}^q (\hat{X}_{li} - X_{li})\frac{\partial f(X_i^*, \boldsymbol{\beta})}{\partial x_l}$ with $X_i^*$ being a point between $X_i$ and $\hat{X}_i$. Condition (A7) and (30) ensure that

$$F_1 = O_P\left((h^4 + n^{-1/2}h^{-1}\log n)^2 \sum_{i=1}^n Y_i^2\right)$$

$$= o_P(n^{1/2}),$$

$$F_2 = O_P\left((h^4 + n^{-1/2}h^{-1}\log n)^2 \left\{\sum_{i=1}^n \left(\sum_{l=1}^q \frac{\partial f(X_i^*, \boldsymbol{\beta})}{\partial x_l}\right)^2\right\}\right)$$

$$= O_P\left((h^4 + n^{-1/2}h^{-1}\log n)^2 \left\{\sum_{i=1}^n C\right\}\right)$$

$$= o_P(n^{1/2}).$$

Arguing as for $F_1$ and $F_2$, we can easily obtain $F_3 = o_P(n)$ and $F_4 = o_P(n^{1/2})$. For $F_5$, given the set $\mathcal{S} = \{(X_i, U_i), 1 \le i \le n\}$, we can derive that

$$
\begin{aligned}
F_5 &= E(F_5|\mathcal{S}) + O_P(\text{Var}(F_5|\mathcal{S})) \\
&= O_P(\text{Var}(F_5|\mathcal{S})) \\
&= O_P\left(\sum_{i=1}^{n}(f(\hat{X}_i, \boldsymbol{\beta}) - f(X_i, \boldsymbol{\beta}))^2\right) \\
&= o_P(n^{1/2}).
\end{aligned}
$$

PROOF OF (ii). Denote the estimating function vector $(R_n^1(\boldsymbol{\beta}), \ldots, R_n^p(\boldsymbol{\beta}))^\tau$ by $R_n(\boldsymbol{\beta})$ and $(G_n^1(\boldsymbol{\beta}), \ldots, G_n^p(\boldsymbol{\beta}))^\tau$ by $G_n(\boldsymbol{\beta})$. From the mean-value theorem, there exists a $\boldsymbol{\beta}^*$ such that

$$
\begin{aligned}
G_n(\boldsymbol{\beta}^0) &= G_n(\hat{\boldsymbol{\beta}}_n) + \frac{\partial G_n(\boldsymbol{\beta}^*)}{\partial \boldsymbol{\beta}}(\boldsymbol{\beta}^0 - \hat{\boldsymbol{\beta}}) \\
&= \frac{\partial G_n(\boldsymbol{\beta}^*)}{\partial \boldsymbol{\beta}}(\boldsymbol{\beta}^0 - \hat{\boldsymbol{\beta}}),
\end{aligned}
\tag{31}
$$

where the $(s, k)$ element of matrix $\frac{\partial G_n(\boldsymbol{\beta}^*)}{\partial \boldsymbol{\beta}}$ is

$$
\begin{aligned}
&\frac{\partial G_n(\boldsymbol{\beta}^*)}{\partial \boldsymbol{\beta}}(s, k) \\
&\quad = -\sum_{i=1}^{n} \frac{\partial f(\hat{X}_i, \boldsymbol{\beta}^*)}{\partial \beta_s} \frac{\partial f(\hat{X}_i, \boldsymbol{\beta}^*)}{\partial \beta_k} + \sum_{i=1}^{n}(\hat{Y}_i - f(\hat{X}_i, \boldsymbol{\beta}^*)) \frac{\partial^2 f(X_i, \boldsymbol{\beta}^*)}{\partial \beta_s \, \partial \beta_k} \\
&\qquad + \sum_{i=1}^{n}(\hat{Y}_i - f(\hat{X}_i, \boldsymbol{\beta}^*)) \sum_{l=1}^{q}(\hat{X}_{li} - X_{li}) \frac{\partial^3 f(X_{li}^*, \boldsymbol{\beta}^*)}{\partial \beta_s \, \partial \beta_k \, \partial x_l}
\end{aligned}
$$

and $\boldsymbol{\beta}^*$ lies between $\boldsymbol{\beta}^0$ and $\hat{\boldsymbol{\beta}}_n$, $X_{li}^*$ between $X_{li}$ and $\hat{X}_{li}$. We want to show that

$$
\sup_{\boldsymbol{\beta}} \left| \frac{1}{n} \frac{\partial G_n(\boldsymbol{\beta})}{\partial \boldsymbol{\beta}}(s, k) + \frac{1}{n} \sum_{i=1}^{n} \frac{\partial f(X_i, \boldsymbol{\beta})}{\partial \beta_s} \frac{\partial f(X_i, \boldsymbol{\beta})}{\partial \beta_k} \right| = o_P(1).
\tag{32}
$$

It suffices to prove that

$$
\sup_{\boldsymbol{\beta}} \left| \frac{1}{n} \sum_{i=1}^{n} \frac{\partial f(\hat{X}_i, \boldsymbol{\beta})}{\partial \beta_s} \frac{\partial f(\hat{X}_i, \boldsymbol{\beta})}{\partial \beta_k} - \frac{1}{n} \sum_{i=1}^{n} \frac{\partial f(X_i, \boldsymbol{\beta})}{\partial \beta_s} \frac{\partial f(X_i, \boldsymbol{\beta})}{\partial \beta_k} \right| = o_P(1),
$$

$$
\sup_{\boldsymbol{\beta}} \left| \frac{1}{n} \sum_{i=1}^{n}(\hat{Y}_i - f(\hat{X}_i, \boldsymbol{\beta})) \frac{\partial^2 f(X_i, \boldsymbol{\beta})}{\partial \beta_s \, \partial \beta_k} \right| = o_P(1),
$$

$$
\sup_{\boldsymbol{\beta}} \left| \frac{1}{n} \sum_{i=1}^{n}(\hat{Y}_i - f(\hat{X}_i, \boldsymbol{\beta})) \sum_{l=1}^{q}(\hat{X}_{li} - X_{li}) \frac{\partial^3 f(X_{li}^*, \boldsymbol{\beta})}{\partial \beta_s \, \partial \beta_k \, \partial x_l} \right| = o_P(1).
$$

Invoking the regularity condition (A1) on $f$ and mimicking the proof of Lemma A.1, the second and third equation can be easily proved. There exist two points $X_i^*$ and $X_i^{**}$, so that the left-hand side of the first equation is bounded by

$$\sup_{\boldsymbol{\beta}} \left| \frac{1}{n} \sum_{i=1}^n \sum_{l=1}^q \frac{\partial^2 f(X_i^*, \boldsymbol{\beta})}{\partial \beta_s \, \partial x_l} \frac{\partial f(X_i, \boldsymbol{\beta})}{\partial \beta_k} (\hat{X}_{li} - X_{li}) \right|$$

$$+ \sup_{\boldsymbol{\beta}} \left| \frac{1}{n} \sum_{i=1}^n \sum_{l=1}^q \frac{\partial f(X_i, \boldsymbol{\beta})}{\partial \beta_s} \frac{\partial^2 f(X_i^{**}, \boldsymbol{\beta})}{\partial \beta_k \, \partial x_l} (\hat{X}_{li} - X_{li}) \right|$$

$$+ \sup_{\boldsymbol{\beta}} \left| \frac{1}{n} \sum_{i=1}^n \left( \sum_{l=1}^q \frac{\partial^2 f(X_i^*, \boldsymbol{\beta})}{\partial \beta_s \, \partial x_l} (\hat{X}_{li} - X_{li}) \right) \right.$$

$$\left. \times \left( \sum_{l=1}^q \frac{\partial^2 f(X_i^{**}, \boldsymbol{\beta})}{\partial \beta_k \, \partial x_l} (\hat{X}_{li} - X_{li}) \right) \right|.$$

Again, by condition (A1) and Lemma A.1, every term above is of order $o_P(1)$. Recalling the definition of $\Lambda$ provided in (11) and combining this with (32) and the consistency of $\hat{\boldsymbol{\beta}}$, we have

$$\frac{1}{n} \frac{\partial G_n(\boldsymbol{\beta}^*)}{\partial \boldsymbol{\beta}} = \Lambda + o_P(1). \tag{33}$$

From Proposition 1, $G_n(\boldsymbol{\beta}) = R_n(\boldsymbol{\beta}) + o_P(n^{1/2})$. Then, (31) and (33) ensure that

$$\begin{aligned}
\sqrt{n}(\hat{\boldsymbol{\beta}} - \boldsymbol{\beta}^0) &= (\Lambda + o_P(1))^{-1} \frac{1}{\sqrt{n}} G_n(\boldsymbol{\beta}^0) \\
&= (\Lambda + o_P(1))^{-1} \left( \frac{1}{\sqrt{n}} R_n(\boldsymbol{\beta}^0) + o_P(1) \right) \\
&= \Lambda^{-1} \frac{1}{\sqrt{n}} R_n(\boldsymbol{\beta}^0) + o_P(1),
\end{aligned} \tag{34}$$

which entails the result. $\quad\square$

PROOF OF THEOREM 2. With $X_0 = 1$ and $\phi_0(u) = 1$, we can easily see, first, that $\breve{\sigma}_k^2 - \sigma_k^2 \geq 0$ is equivalent to

$$\sigma^2 (\Lambda^{-1})_{kk} + \frac{\sigma^2}{\mathbb{E}[Y]} \beta_k I(k=0) + \frac{1}{4} \beta_k^2 \mathbb{E} \left( \frac{Y - \mathbb{E}[Y]}{\mathbb{E}[Y]} - \frac{X_k - \mathbb{E}[X_k]}{\mathbb{E}[X_k]} \right)^2$$

$$+ \beta_k^2 \mathbb{E} \left( \frac{\widetilde{Y} - Y}{\mathbb{E}[Y]} - \frac{\widetilde{X}_k - X_k}{\mathbb{E}[X_k]} \right)^2$$

$$\leq \sigma^2 (\Lambda^{-1})_{kk} + \sigma^2 (\Lambda^{-1})_{kk} \operatorname{Var}[\psi(U)] + \beta_k^2 \frac{\mathbb{E}[X_k^2]}{(\mathbb{E}[X_k])^2} \operatorname{Var}[\psi(U) - \phi_k(U)].$$

After some algebraic calculations, it is also equivalent to

$$\frac{\sigma^2}{\mathbb{E}[Y]}\beta_k I(k=0) + \frac{1}{4}\beta_k^2 \mathbb{E}\left(\frac{Y}{\mathbb{E}[Y]} - \frac{X_k}{\mathbb{E}[X_k]}\right)^2$$

$$+ \beta_k^2 \mathbb{E}\left(\frac{Y}{\mathbb{E}[Y]} - \frac{X_k}{\mathbb{E}[X_k]}\right)^2 \mathrm{Var}[\psi(U)]$$

$$+ 2\beta_k^2 \mathbb{E}\left(\left(\frac{Y}{\mathbb{E}[Y]} - \frac{X_k}{\mathbb{E}[X_k]}\right)\frac{X_k}{\mathbb{E}[X_k]}\right)\mathbb{E}[\psi^2(U) - \psi(U)\phi_k(U)]$$

$$\leq \sigma^2 (\Lambda^{-1})_{kk} \mathrm{Var}[\psi(U)].$$

Then, multiplying the two sides of the above inequality by $(\mathbb{E}[Y])^2$, and noticing that $\mathbb{E}[Y^2] = \boldsymbol{\beta}^\tau \Lambda \boldsymbol{\beta} + \sigma^2$, $(\mathbb{E}[Y])^2 = \boldsymbol{\beta}^\tau \Lambda_0 \Lambda_0^\tau \boldsymbol{\beta}$, and $\mathbb{E}[YX_k]\mathbb{E}[Y] = \boldsymbol{\beta}^\tau(\frac{1}{2}\Lambda_k \times \Lambda_0^\tau + \frac{1}{2}\Lambda_0\Lambda_k^\tau)\boldsymbol{\beta}$, we can obtain that $\breve{\sigma}_k^2 - \sigma_k^2 \geq 0$ is equivalent to

$$(35) \qquad \boldsymbol{\beta}^\tau \mathbf{C}_k \boldsymbol{\beta} \leq 0, \qquad 0 \leq k \leq p,$$

where $\mathbf{C}_k$ is given as

$$\mathbf{C}_k = \left\{(-3 - 4\mathbb{E}[\psi^2(U)] + 8\mathbb{E}[\psi(U)\phi_k(U)])\frac{\mathbb{E}[X_k^2]}{(\mathbb{E}[X_k])^2}\beta_k^2\right.$$

$$\left. - 4\sigma^2 \mathrm{Var}[\psi(U)](\Lambda^{-1})_{kk}\right\}\Lambda_0\Lambda_0^\tau$$

$$+ \sigma^2(-3 + 4\mathbb{E}[\psi^2(U)])e_k e_k^\tau + 2\sigma^2(e_k\Lambda_k^\tau + \Lambda_k e_k^\tau)I(k=0)$$

$$+ (3 - 4\mathbb{E}[\psi(U)\phi_k(U)])\frac{1}{\mathbb{E}[X_k]}\beta_k^2(\Lambda_k\Lambda_0^\tau + \Lambda_0\Lambda_k^\tau)$$

$$+ (-3 + 4\mathbb{E}[\psi^2(U)])\beta_k^2 \Lambda.$$

Clearly, the matrices $\mathbf{C}_k$s are symmetric. In the following, we verify that the set about $\boldsymbol{\beta}$ satisfying the above $p+1$ inequalities is not empty. Because $\mathbf{C}_k$ involves the parameter $\beta_k$, (35) is, indeed, not a standard quadric form. Instead, we can investigate the equivalent standard one, which is

$$\boldsymbol{\beta}_k^\tau \mathbf{C}_k \boldsymbol{\beta}_k \leq 0, \qquad 0 \leq k \leq p$$

with $\boldsymbol{\beta}_k = \boldsymbol{\beta}/\beta_k$. $\boldsymbol{\beta}_k$ varies with $k$, but their directions are the same. This means that it suffices to consider (35). To prove (35) nonempty, we verify that there exists at least one negative eigenvalue of matrices $\mathbf{C}_k$s. For any $n \times n$ matrix, denote by $\lambda_{\max}(\cdot)$ and $\lambda_{\min}(\cdot)$ its largest and smallest eigenvalues, respectively. To prove this, we only need to prove that

$$(36) \qquad \lambda_{\min}(\mathbf{C}_k) \leq 0, \qquad 0 \leq k \leq p.$$

Note that, for any given symmetric matrix $B$, it always holds that $\lambda_{\min}(B) \leq B_{ii} \leq \lambda_{\max}(B)$. So, we obtain that,

$$
\begin{aligned}
\lambda_{\min}(\mathbf{C}_k) &\leq (\mathbf{C}_k)_{00} \\
&= \beta_k^2 \left(1 - \frac{\mathbb{E}[X_k^2]}{(\mathbb{E}[X_k])^2}\right) (4\mathbb{E}[\psi^2(U)] - 8\mathbb{E}[\psi(U)\phi_k(U)] + 3) \\
&\quad - 4\sigma^2(\Lambda^{-1})_{kk} \operatorname{Var}[\psi(U)] + \sigma^2(4\operatorname{Var}[\psi(U)] + 5)I(k = 0).
\end{aligned}
$$

If the distorting functions $\psi(u)$ and $\phi_k(u)$ satisfy

$$
4\mathbb{E}[\psi^2(U)] - 8\mathbb{E}[\psi(U)\phi_k(U)] + 3 \geq 0
$$

and

$$
\operatorname{Var}[\psi(U)] \geq \frac{5}{4((\Lambda^{-1})_{00} - 1)},
$$

then $(\mathbf{C}_k)_{00} \leq 0$, and then $\lambda_{\min}(\mathbf{C}_k) \leq 0$. (36) is proved; that is, the set confined by the conditions of Theorem 2 is not empty. $\square$

PROOF OF THEOREM 3.  By some elementary calculations, we can have

$$
\begin{aligned}
G_{n,i}(\boldsymbol{\beta}) &= (Y_i - f(X_i, \boldsymbol{\beta})) \frac{\partial f(X_i, \boldsymbol{\beta})}{\partial \beta} \\
&\quad - (f(\hat{X}_i, \boldsymbol{\beta}) - f(X_i, \boldsymbol{\beta})) \frac{\partial f(X_i, \boldsymbol{\beta})}{\partial \beta} \\
&\quad + (\hat{Y}_i - Y_i) \frac{\partial f(X_i, \boldsymbol{\beta})}{\partial \beta} \\
&\quad + (\hat{Y}_i - f(\hat{X}_i, \boldsymbol{\beta})) \left(\frac{\partial f(\hat{X}_i, \boldsymbol{\beta})}{\partial \beta} - \frac{\partial f(X_i, \boldsymbol{\beta})}{\partial \beta}\right) \\
&\equiv (Y_i - f(X_i, \boldsymbol{\beta})) \frac{\partial f(X_i, \boldsymbol{\beta})}{\partial \beta} - E_{i1} + E_{i2} + E_{i3}.
\end{aligned}
$$
(37)

We now show that

$$
\max_{1 \leq i \leq n} |G_{n,i}(\boldsymbol{\beta})| = o_P(n^{1/2}).
$$
(38)

It is well known that, for any sequence of independent and identically distributed random variables $\{Z_i, 1 \leq i \leq n\}$ with $\mathbb{E}[Z_i^2] < \infty$,

$$
\max_{1 \leq i \leq n} |Z_i|/\sqrt{n} \to 0, \qquad \text{almost surely.}
$$

This implies that $\max_{1 \leq i \leq n} |(Y_i - f(X_i, \boldsymbol{\beta})) \frac{\partial f(X_i, \boldsymbol{\beta})}{\partial \beta}| = o_P(n^{1/2})$, because $\mathbb{E}[(Y_i - f(X_i, \boldsymbol{\beta})) \frac{\partial f(X_i, \boldsymbol{\beta})}{\partial \beta}]^2 < \infty$. Using Taylor expansion to the quantity $f(\hat{X}_i, \boldsymbol{\beta}) -$

$f(X_i, \boldsymbol{\beta})$, with respect to $X$ and by condition (A1), $\max_{1 \le i \le n} |E_{i1}| = o_P(n^{1/2})$ follows from

$$
\begin{aligned}
|E_{i1}| &\le \sum_{l=1}^{q} \left| \frac{\partial f(X_i^*, \boldsymbol{\beta})}{\partial x_l} (\hat{X}_{li} - X_{li}) \frac{\partial f(X_i, \boldsymbol{\beta})}{\partial \beta} \right| \\
&\le C \sum_{l=1}^{q} \left| \left( \frac{g_l(U_i) \hat{p}(U_i)}{\hat{g}_l(U_i) p(U_i)} \frac{\bar{\bar{X}}_l}{\mathbb{E}[X_l]} - 1 \right) X_{li} \frac{\partial f(X_i, \boldsymbol{\beta})}{\partial \beta} \right| \\
&= C \sum_{l=1}^{q} \left| \left( \frac{g_l(U_i) \hat{p}(U_i)}{\hat{g}_l(U_i) p(U_i)} - 1 \right) \right. \\
&\qquad + \left. \left( \frac{g_l(U_i) \hat{p}(U_i)}{\hat{g}_l(U_i) p(U_i)} - 1 \right) \frac{\bar{\bar{X}}_l - \mathbb{E}[X_l]}{\mathbb{E}[X_l]} + \frac{\bar{\bar{X}}_l - \mathbb{E}[X_l]}{\mathbb{E}[X_l]} \right| \left| X_{li} \frac{\partial f(X_i, \boldsymbol{\beta})}{\partial \beta} \right| \\
&\le C \sum_{l=1}^{q} \left| \sup_u \left| \frac{g_l(u) \hat{p}(u)}{\hat{g}_l(u) p(u)} - 1 \right| \right. \\
&\qquad + \left. O_P(n^{-1/2}) \sup_u \left| \frac{g_l(u) \hat{p}(u)}{\hat{g}_l(u) p(u)} - 1 \right| + O_P(n^{-1/2}) \right| \left| X_{li} \frac{\partial f(X_i, \boldsymbol{\beta})}{\partial \beta} \right| \\
&= \sum_{l=1}^{q} \{ O_P(h^4 + n^{-1/2} h^{-1} \log n) + O_P(n^{-1/2}) \} \left| X_{li} \frac{\partial f(X_i, \boldsymbol{\beta})}{\partial \beta} \right|.
\end{aligned}
$$

Applying similar techniques to those used in the analysis of $E_{i1}$, we have $\max_{1 \le i \le n} |E_{i2}| = o_P(n^{1/2})$ and $\max_{1 \le i \le n} |E_{i3}| = o_P(n^{1/2})$.

Further, applying the results of Proposition 1 and following the same arguments as were used in the proof of expression (2.14) in Owen (1990), we have

$$
\lambda = O_P(n^{-1/2}). \tag{39}
$$

Applying Taylor expansion to (16), and invoking (38), (39) and (10) of Proposition 1, we can verify that

$$
\hat{l}(\boldsymbol{\beta}) = 2 \sum_{i=1}^{n} \left( \lambda^\tau G_{n,i}(\boldsymbol{\beta}) - \frac{1}{2} \{ \lambda^\tau G_{n,i}(\boldsymbol{\beta}) \}^2 \right) + o_P(1). \tag{40}
$$

By (17), it follows that

$$
\begin{aligned}
0 &= \frac{1}{n} \sum_{i=1}^{n} \frac{G_{n,i}(\boldsymbol{\beta})}{1 + \lambda^\tau G_{n,i}(\boldsymbol{\beta})} \\
&= \frac{1}{n} \sum_{i=1}^{n} G_{n,i}(\boldsymbol{\beta}) - \frac{1}{n} \sum_{i=1}^{n} G_{n,i}(\boldsymbol{\beta}) G_{n,i}^\tau(\boldsymbol{\beta}) \lambda + \frac{1}{n} \sum_{i=1}^{n} \frac{G_{n,i}(\boldsymbol{\beta}) \{ \lambda^\tau G_{n,i}(\boldsymbol{\beta}) \}^2}{1 + \lambda^\tau G_{n,i}(\boldsymbol{\beta})}.
\end{aligned}
$$

The application of (38), (39) and (10) of Proposition 1 yields

$$
(41) \qquad \lambda = \left( \sum_{i=1}^{n} G_{n,i}(\boldsymbol{\beta}) G_{n,i}^{\tau}(\boldsymbol{\beta}) \right)^{-1} \sum_{i=1}^{n} G_{n,i}(\boldsymbol{\beta}) + o_P(n^{-1/2}).
$$

Equation (40), together with (41), obtains the decomposition

$$
\begin{aligned}
(42) \qquad \hat{l}(\boldsymbol{\beta}) = {} & n \left( \frac{1}{n} \sum_{i=1}^{n} G_{n,i}(\boldsymbol{\beta}) \right)^{\tau} \left( \frac{1}{n} \sum_{i=1}^{n} G_{n,i}(\boldsymbol{\beta}) G_{n,i}^{\tau}(\boldsymbol{\beta}) \right)^{-1} \\
& \times \left( \frac{1}{n} \sum_{i=1}^{n} G_{n,i}(\boldsymbol{\beta}) \right) + o_P(1).
\end{aligned}
$$

From Proposition 1 and (42), we have

$$
(43) \qquad \hat{l}(\boldsymbol{\beta}) = R_n^{\tau}(\boldsymbol{\beta})(R_n(\boldsymbol{\beta}) R_n^{\tau}(\boldsymbol{\beta}))^{-1} R_n(\boldsymbol{\beta}) + o_P(1).
$$

Since $R_n(\boldsymbol{\beta})$ is a sum of independent and identically distributed random variables, the result of Theorem 3 is immediately achieved by central limit theorems.  $\square$

**Acknowledgments.**  The authors thank an Associate Editor and two referees for their constructive comments and suggestions, which led to great improvements over an early manuscript.

X. CUI
L. LIN
SCHOOL OF MATHEMATICS
SHANDONG UNIVERSITY
JINAN, SHANDONG PROVINCE, 250100
P.R. CHINA
E-MAIL: xcui@mail.sdu.edu.cn
        linlu@sdu.edu.cn

W. GUO
CENTER FOR CLINICAL EPIDEMIOLOGY
    AND BIOSTATISTICS
UNIVERSITY OF PENNSYLVANIA
    SCHOOL OF MEDICINE
613 BLOCKLEY HALL
423 GUARDIAN DRIVE
PHILADELPHIA, PENNSYLVANIA 19104-6021
USA
E-MAIL: wguo@mail.med.upenn.edu

L. ZHU
FSC1207
FONG SHU CHUEN BUILDING
DEPARTMENT OF MATHEMATICS
THE HONG KONG BAPTIST UNIVERSITY
KOWLOON TONG
HONG KONG
P.R. CHINA
E-MAIL: lzhu@hkbu.edu.hk